\documentclass[default,a4paper,twoside,onecolumn,10pt]{sn-jnl}

\jyear{2022}%

\theoremstyle{thmstyleone}%
\theoremstyle{thmstyletwo}%

\theoremstyle{thmstylethree}%
\usepackage{usrextra} 
\usepackage{lscape}
\usepackage{psfrag}
\usepackage{lscape}
\usepackage{color}  
\usepackage{fancyhdr}  
\usepackage{natbib}
\usepackage{hd}
\usepackage{mathtools, cuted}
\usepackage{lipsum}
\begin{document}

\title[ ]{An \textsl{in silico}-based review on anisotropic hyperelastic constitutive models for soft biological tissues}{}

\author[1]{H\"usn\"u Dal\footnote{Corresponding author: dal@metu.edu.tr}}

\author[1]{Alp Ka\u{g}an A\c{c}an}\email{alpacan@metu.edu.tr}

\author[2]{Ciara Durcan}\email{998131@swansea.ac.uk}

\author[2]{Mokarram Hossain }\email{mokarram.hossain@swansea.ac.uk}

%
%
%

\affil[1]{\orgdiv{Department of Mechanical Engineering}, \orgname{Middle East Technical University}, \orgaddress{\state{Dumlup\i nar Bulvar\i~1, 06800 \c Cankaya, Ankara }, \country{T\"urkiye}}}

\affil[2]{\orgdiv{Zienkiewicz Centre for Computational Engineering, College of Engineering}, \orgname{Swansea University}, \orgaddress{\state{Swansea, SA1 8EN}, \country{United Kingdom}}}


\abstract{We review nine invariant and dispersion-type anisotropic hyperelastic constitutive models for soft biological tissues based on their fitting performance to experimental data from three different human tissues. For this, we used a hybrid multi-objective optimization procedure. A genetic algorithm is devised to generate the initial guesses followed by a gradient-based search algorithm. The constitutive models are then fit to a set of uniaxial and biaxial tension experiments conducted on tissues with differing fiber orientations. For the \textsl{in silico} investigation, experiments conducted on aneurysmatic abdominal aorta, linea alba, and rectus sheath tissues are utilized. Accordingly, the models are ranked with respect to an objective normalized quality of fit metric. Finally, a detailed discussion is carried out on the fitting performance of each model. This work provides a valuable quantitative comparison of various anisotropic hyperelastic models, the findings of which can aid those modeling the behavior of soft tissues in selecting the best constitutive model for their particular application.}




\maketitle

\section{Introduction}
\label{chp:b1} 
The physiological functions of many soft biological tissues are related to their material composition and mechanical properties. Examples of these functions include blood flow through the heart and blood vessels~\cite{peskin1977,mcculloch2020,carvalho2021}, peristalsis of the gastrointestinal (GI) tract~\cite{lew1971,alokaily2019}, and even the transmission of sound within the ear cavity~\cite{zhang2020}. The mechanical properties of tissues can also be a predictor of their pathophysiology, such as the stiffness of arteries being a successful indicator of cardiovascular disease, or the elastic and viscoelastic parameters of skin being linked to the occurrence of autoimmune disorders~\cite{kuwahara2008,kalra2016}. To successfully model the physiological functions of soft tissues and to investigate any pathophysiology, the characteristics of the material must first be established. Constitutive modeling can be used to analytically capture the mechanical behavior of these biological tissues.  Such an ansatz forms the basis for developing computational models which consider the geometry of the tissue or organ to study how shape contributes to its function. Baseline hyperelasticity is a precursor in biomechanics to the modeling of more advanced phenomena such as viscoelasticity~\cite{cansiz+dal+etal15}, electroelasticity~\cite{dal+goktepe+etal13} and electro-viscoelasticity~\cite{cansiz+dal+etal17,cansiz+dal+etal18}, to mention a few. In addition to this, computational modeling can be used to study fluid-structure interaction scenarios, such as urine traveling down the urinary tract~\cite{zheng2021}, or air filling up the alveoli of the lungs~\cite{wall2008}. However, before this is possible, the behavior of the solid structures must be determined and successfully simulated. A variety of applications exist for the computational models mentioned, including within the medical device industry, where  the mathematical models can be used to aid in the design of devices that interact with the human body, such as stents and endoscopy devices~\cite{finotello2021,tian2021}. The use of the models in a finite element implementation can give greater insight into the devices' mechanisms, highlight any potential risks for the patient, and conserve a range of resources such as time and manufacturing materials. Further to this, computational models of soft tissues can be used for surgical simulations; the training technique which allows surgeons to practice a procedure an infinite number of times without risking the health or lives of patients~\cite{misra2008,goulette2015,afshari2017}.

It has previously been well established that soft tissues respond elastically to a given force and past a certain point, i.e. under large strains, behave non-linearly~\cite{fung1967,beatty1987,weiss1996}. For this, and when the time-dependent properties of the tissue are not considered, hyperelastic modeling, similar to that for polymers, is required~\cite{steinmann2012,hossain2012}. For polymers, isotropic, hyperelastic models are commonly used as these substances often have a homogeneous material composition, unless specifically designed otherwise. Biological tissues, however, are generally observed through experimentation to be anisotropic. For fibrous tissues, this is often associated with different densities and/or orientations of the tissue's main constituent fibers, resulting in a differing response of the tissue depending on its loading direction. The anisotropy of biological tissues, hard or soft, can normally be linked with the physiological function of the organ when considering the mechanical characteristics of each different fiber.
There are many types of tests to characterize the anisotropic hyperelastic behavior of biological tissues, and amongst the most common include uniaxial tension tests in two or more directions~\cite{myers2010}, biaxial tension tests which simultaneously load in two directions~\cite{geest2006,huntington2019}, pure shear tension tests in two directions~\cite{masri2018}, inflation-extension tests for tubular tissues~\cite{sommer2010} and compression tests in different directions~\cite{bol2014}. Experiments conducted on soft biological tissues can be difficult due to their delicate nature in addition to practical limitations including the need for highly sensitive data acquisition equipment to capture their comparatively low internal forces. This makes the experimental studies conducted on soft tissues extremely valuable in the field of biomechanics. 
Such studies conducted using animal soft tissues are more prevalent than with human tissue due to the former being more readily available and having fewer associated ethical constraints. To date, some non-exhaustive examples of the experimental data available on animal soft tissues includes the porcine aorta~\cite{pena2015}, rodent and porcine brain~\cite{macmanus2017}, rodent vagina~\cite{huntington2019}, porcine liver~\cite{mattei2014}, ovine oesophagi~\cite{sommer2013}, swine tongue~\cite{wang2018a}, equine urethra~\cite{natali2016}, porcine skeletal muscle~\cite{takaza2013}, rabbit skin~\cite{sherman2017} and porcine stomach~\cite{jia2015}. 

Animal tissue may provide a good representation of how human tissue responds to certain loads or strains, but the data from such experiments cannot be used to accurately model human tissues, especially for use within medical applications. Despite being less prevalent than animal studies, there are also a wide variety of experimental studies on human soft tissues, with investigations being conducted both {\it in vivo} and {\it ex vivo}. In this review, we will focus purely on human data from {\it ex vivo} experiments as these allow for large  deformations that clearly display the non-linear elastic behavior of the tissues. Experimental data for healthy human soft tissues can be found for the brain~\cite{budday2017}, spleen~\cite{kemper2012}, skin~\cite{groves2013}, aorta~\cite{ferrara2016,vignali2021}, bile duct~\cite{girard2019}, oesophagus ~\cite{durcan2022a,durcan2022b}, GI tract~\cite{egorov2002}, liver~\cite{untaroiu2015}, cervical tissue~\cite{myers2010}, while also for diseased tissue such as aneurysmatic aortas~\cite{sommer2016} and diabetic foot plantar tissues~\cite{pai2010}. 
As mentioned previously, and as can be witnessed from the time-independent stress-strain relationships of these experimental studies, at large strains soft biological tissues exhibit a non-linear hyperelastic material response, therefore their behavior cannot be simply captured by two constants of linear elasticity. A strain energy density function (SEDF) must be used to model the hyperelastic response of a soft tissue to certain loads or deformations, and  SEDFs exist in many different forms. There are a variety of approaches including the very popular classical models by Fung et al.~\cite{Fung1979} and Holzapfel et al.~\cite{holzapfel2000}, as well as more recently developed models~\cite{chaimoon2019,chanda2020}, and those which describe certain time-independent phenomena such as the Mullins effect~\cite{rebouah2014}. Micromechanical models provide insight into the mechanical role of the different components that comprise a material. In these models, every parameter is linked to a physical phenomenon, with the aim to predict the behavior of the tissue based purely on its composition. Micromechanical dispersion models give further insight into the orientation and distribution of the tissues' main constituent fibers. The parameter in these models associated with the dispersion of the fibers should ideally be identified using microstructural data obtained from the tissue in question~\cite{Holzapfel2015}. Phenomenological models, on the other hand, utilize material parameters that do not have a physical meaning, these, however, are often found to capture the experimental behavior very well and may be more relevant for a specific application or loading domain~\cite{Fung1979}. 

A number of comprehensive reviews have been carried out on the topic of hyperelastic constitutive models. Steinmann et al.~\cite{steinmann2012}, Hossain et al.~\cite{hossain2013, hossain2015} and Dal et al.~\cite{Dal2021} reviewed the range of hyperelastic models in literature. A review of the models used to simulate the behavior of anisotropic hyperelastic soft biological tissues was conducted by Chagnon et al.~\cite{chagnon2015}. This review presents an extensive and detailed list of many prominent models used in the field.  Wex et al.~\cite{wex2015} provided a review of the most commonly used isotropic constitutive models for biological tissues. Mihai et al.~\cite{mihai2015} reviewed models used for the brain and fat tissues, and successfully compared their performance ability under different loading conditions. Furthermore, Rosen and Jiang~\cite{rosen2019} investigated current isotropic hyperelastic laws able to model tissue-mimicking materials characterized through shear wave elastography. Martins et al.~\cite{martins2006} carried out a concise review on isotropic models subjected to uniaxial loading. Wang et al.~\cite{wang2018b} added a level of complexity and reviewed the state of modeling the compressibility of soft tissues. Holzapfel et al.~\cite{holzapfel2019} looked at purely dispersion-type models. Bhattarai et al.~\cite{bhattarai2021} reviewed the modeling of only one tissue, the large intestine. Recently, Cheng and Zhang~\cite{cheng2018} provided a framework for the derivation of hyperelastic models, both isotropic and anisotropic, in the context of biomaterials.

The aforementioned reviews exclude a one-to-one quantitative comparison of models used for soft tissues. In this work, we provide a derivation of invariant- and dispersion-type models, an outline of the most prominent constitutive models in literature, and an \textsl{in silico} comparison of the models' ability to simulate the anisotropic hyperelastic behavior of a variety of human tissues. The data for this is obtained from experimental studies of aortic aneurysms and different regions of the healthy abdominal wall. This work provides a detailed framework of how to utilize the constitutive models, while also providing an in-depth critical analysis of each models' performance ability according to different data sets and loading domains. The findings of this review can be used in research or industry to support and inform the selection of the best constitutive model for a specific soft tissue application.  

We pursue our investigation based on the  datasets proposed by Niestrawska et al.~\cite{niestrawska2016}, Cooney et al.~\cite{cooney2016}, and Martins et al.~\cite{Martins2012}  for abdominal aortic aneurysm (AAA) tissue, the linea alba and the rectus sheath, respectively. AAA is a pathological condition of the abdominal aorta which results in local bulging, while the linea alba and rectus sheath are sections of the abdominal wall. One reason for selecting these tissues is the fact that they are often subjected to surgical treatment, therefore an accurate representation of their mechanical response leads to representative numerical analysis for surgical procedures and the design of proper medical devices. In addition, each tissue presents a unique stress-stretch response, thus the choice of these tissues allows the models to be assessed across different mechanical behaviors. The overall output of the paper is twofold. Firstly, the applicability of each constitutive model to a particular type of tissue will be assessed, specifically the aptness of various anisotropic fiber distribution functions to the three tissues mentioned above. Secondly, the fitting performance of various dispersion-type anisotropic constitutive models will be compared. In this study, we will
\begin{itemize}
 \item  outline nine anisotropic hyperelastic models and derive the closed form stress expressions for the uniaxial tension (UT) and equibiaxial tension (ET) deformation states based on the classical incompressibility assumption,
\item implement the aforementioned  constitutive models into a hybrid genetic-gradient search algorithm outlined in Dal et al.~\cite{Dal2021},
\item identify the material parameters of each model (i) with respect to an ET dataset for AAA tissue \cite{niestrawska2016}, (ii) a UT dataset for the linea alba~\cite{cooney2016}, and (iii) a UT dataset for the rectus sheath~\cite{Martins2012},
\item sort all the models with respect to an objective \textsl{quality of fit} metric according to their fitting performance to the AAA tissue dataset \cite{niestrawska2016}, the linea alba dataset \cite{cooney2016}, and the rectus sheath dataset \cite{Martins2012}.
\item The identified parameters and the quality of fit values of each constitutive model will be presented in tabular form. Also, the stress-strain curves for each constitutive model will be provided in separate graphs.
\end{itemize}
The paper is organized as follows: the mathematical preliminaries for the kinematics of incompressible anisotropic hyperelastic solids are presented in Section \ref{chp-b2}. In Section \ref{chp:b3}, nine anisotropic hyperelastic constitutive models for soft biological tissues are summarized. Section \ref{chp:b4} outlines the experimental procedures for the characterization of the quasi-static hyperelastic response of the soft human tissues selected in this study. In Section \ref{chp:b5}, the parameter identification  procedure and  quality of fit metric are provided. Finally, a detailed discussion of the results is presented in Section \ref{chp:b6}.


\section{Fundamentals of hyperelastic deformable solids} 
\label{chp-b2}
In this section, the kinematics, state variables and stress expressions for an anisotropic hyperelastic continuum will be introduced. The mathematical framework of invariant and dispersion-type anisotropic hyperelastic formulations will be briefly discussed.  
\subsection{Geometric mappings and the field variables}
\subsubsection{Kinematics}
{\setlength{\parindent}{0cm}
Let the deformation map $\Bvarphi(\BX,t)$ represent the motion of a deformable solid body.   It  maps the reference/Lagrangian configuration $\BX\in\SB_0$ of material point onto the current/Eulerian configuration of material points $\Bx=\Bvarphi_t(\BX)$ at time $t\in \ST \subset \IR_+$. The \textsl{deformation gradient}
\begin{equation}
 \BF:T_X\SB_0\rightarrow T_x\SB ;\quad \BF := \pp{\phit(\BX)}{\BX}
\end{equation}
maps the unit tangent of the reference configuration onto the spatial counterpart in the current configuration where $T_X\SB_0$ and $T_x\SB$ denote the tangent spaces in the reference and current configurations, respectively. Additionally, the co-tangent spaces in the reference and current manifolds are represented as $\CSBX$ and $\CSSx$, respectively. The normal map  between the unit normals in the undeformed and deformed configurations is defined as
\begin{equation}\label{unit_normal}
 \BF^{-T}:\CSBX\rightarrow \CSSx ;\quad \Bn=\BF^{-T}\BN \,,
\end{equation}
where $\Bn$ is the unit normal vector in the reference configuration and $\BN$ is the unit normal vector in the current configuration. In this sequence, let line, area, and volume elements in the Lagrangian configurations to be denoted as $d\BX$, $d\BA$, $dV$, respectively. The Eulerian counterparts of these elements are obtained through the deformation gradient $\BF$, its cofactor $\cof [\BF]= \det [\BF]\BF^{-T}$ and its Jacobian $J:=\det [\BF]$
\eb \label{three-mappings}
d\Bx = \BF d\BX~, \qquad d\Ba = \cof [\BF] d\BA~, \qquad dv= \det[\BF] dV\,,
\ee
\inputfig{FIGs/threeFundamentalMaps}{threefundmaps}
see also Figure \ref{three-fund-maps}.  $J:=\det [\BF]>0$ guarantees the nonpenetrable deformations $\Bvarphi$. To define the mappings between the co- and contravariant objects in the Lagrangian and Eulerian manifolds, the undeformed $\SB_0$ and deformed $\SB$ configurations are locally equipped with covariant $\BG$ and current $\Bg$ metric tensors in the neighborhoods $\eulN_{X}$ of $\BX$ and $\eulN_{x}$ of $\Bx$, respectively. Here, the \textsl{right Cauchy-Green tensor} and  the inverse of the \textsl{left Cauchy-Green tensors} 
\eb
\BC = \BF^T\Bg\BF
\qquad \text{and} \qquad
\Bc = \BF^{-T} \BG \BF^{-1} =:\Bb^{-1}
\ee
are defined as the pull back of the current metric $\Bg$ and push-forward of the Lagrangian metric $\BG$, respectively. The inverse of the \textsl{left Cauchy-Green tensor} is commonly known as the \textsl{Finger tensor} and represented as $\Bb=\Bc^{-1}$.
\subsubsection{Principal stretches and invariants}
The principle of material objectivity and the principle of frame indifference require that the energy stored in the hyperelastic material be a function of either principal stretches or invariants of the right Cauchy-Green tensor. The spectral decomposition of right Cauchy-Green tensor and its cofactor reads 
\eb
\BC := \sum_{a=1}^3 \lambda_a^2\BN^a\otimes\BN^a 
\quad \text{and} \quad
\cof[\BC] := \sum_{a=1}^3 \nu_a^2\BN^a\otimes\BN^a \,,
\ee
where 
\eb
\nu_i=J/\lambda_i \quad \text{with}\quad \nu_1=\lambda_2\lambda_3\,, \quad \nu_2=\lambda_3\lambda_1\,, \quad \nu_3=\lambda_1\lambda_2\,, 
\ee
are the \textsl{principal areal stretches}, see Figure \ref{isoinv}. 
\inputfig{FIGs/isoinv}{isoinv}
Moreover, the three isotropic invariants of the right Cauchy-Green tensor are
\eb\label{inv-isotropic}
I_1:= \tr[\BC] , \quad   I_2:=\tr[\cof[\BC]]  ,\quad \textrm{and}  \quad I_3:=\det[\BC]\,.
\ee
The principal stretches and the invariants of the right Cauchy-Green stretch tensor are related as 
\eb\label{stretch-isotropic}
I_1=\lambda_1^2+\lambda_2^2+\lambda_3^2, \quad I_2=\nu_1^2+\nu_2^2 + \nu_3^2, \quad  I_3=J^2=\lambda_1^2\lambda_2^2\lambda_3^2\,.
\ee
For an infinitesimal cubic element, the three isotropic invariants of the right Cauchy-Green stretch tensor are associated with \textsl{linear}, \textsl{areal}, and \textsl{volumetric} stretches in the principal directions. Let the Lagrangian unit vectors $\BM_{1}$ and  $\BM_{2}$ denote the orientation directions of two fiber families for the description of the anisotropic behavior of soft tissues. Four additional invariants that take into account the contribution of the preferential initial fiber alignment are defined as
\begin{equation}\begin{array}{rll}
I_4:= \BM_{1}\cdot\BC\BM_{1}  \qquad \textrm{and}  \qquad I_5:=\BM_{1}\cdot\BC^2\BM_{1} \ , \\[1Ex]
I_6:= \BM_{2}\cdot\BC\BM_{2}  \qquad \textrm{and}  \qquad I_7:=\BM_{2}\cdot\BC^2\BM_{2} \ ,
\label{inv-anisotropic}
\end{array}\end{equation}
where $\lambda^f_{1}=\sqrt{I_4}$ and $\lambda^f_{2}=\sqrt{I_6}$ represent the stretches of the fibers in the orientation directions $\BM_{1}$ and $\BM_{2}$, respectively. Here, the unit orientation vectors of two families of fibers, symmetrically disposed in the $\Be_1-\Be_2$ plane in current configuration are 
\begin{equation}
\Bm_1=\BF\BM_{1} \qquad \textrm{and} \qquad \Bm_2=\BF\BM_{2} \ .
\end{equation}
Moreover, the second order structure tensors can be defined as
\begin{equation}
\BA_1=\BM_{1} \otimes \BM_{1} \qquad \textrm{and} \qquad \BA_2=\BM_{2} \otimes \BM_{2} \ .
\end{equation}
\subsection{Free energy function and the stress expressions}
Hyperelastic materials  are governed by a potential function, the Helmholtz free energy function that describes the energy stored as a result of mechanical deformation. Polymeric materials and soft biological tissues exhibit a distinct response to bulk deformation and shear-type deformations. In this context, it is common practice to split the deformation into  dilatational and volume-preserving parts  
\eb\label{mult-split}
\BF=\BF_{\text{vol}}\bar\BF \,\quad \text{with}\,\quad \BF_{\text{vol}}:=J^{1/3}\Bnone\,.
\ee
The corresponding deformation measure reads
\eb\label{mult-split-C}
\BC=(J^{2/3}\Bnone)\bar\BC \,\quad \text{with}\,\quad \bar\BC = \bar\BF^{T}\bar\BF\,.
\ee
Based on (\ref{mult-split}) and (\ref{mult-split-C}), the free energy function is expressed as 
\begin{equation}\label{split1}
\Psi(\BF,\BA_{1},\BA_{2}) = U(J) + \bar\Psi(\bar\BF,\BA_{1},\BA_{2})
\end{equation}
where $U(J)$ and $\bar\Psi(\bar\BF,\BA_{1},\BA_{2})$ represent the volumetric and isochoric response of the material, respectively.
Further decomposition of the isochoric free energy function was suggested by Holzapfel and Weizsacker~\cite{Holzapfel1998}. They proposed decomposition of the free energy function into isotropic and anisotropic parts,
\begin{equation}\label{split2}
\bar\Psi(\bar\BF,\BA_{1},\BA_{2}) = \bar\Psi_{\text{iso}}(\bar\BF) + \bar\Psi_{\text{ani}}(\bar\BF,\BA_{1},\BA_{2})\,.
\end{equation}
The wavy collagen fibers do not store energy under contractile deformations. Hence, the isotropic ground matrix is active  $\bar\Psi_{\text{iso}}$ under small deformations. On the other hand, the collagen fibers dominate the overall behavior at high stretch levels governed by $\bar\Psi_{\text{ani}}$~\cite{holzapfel2000}.
Following the arguments of G\"ultekin et al.~\cite{Gueltekin2019,Gueltekin2020}, we employ the specific form of the free energy function
\begin{equation}\label{unsplit1}
\Psi(\BF,\BA_{1},\BA_{2}) = U(J) + \tilde\Psi(\BF,\BA_{1},\BA_{2})
\end{equation}
leading to
\begin{equation}\label{unsplit2}
\tilde\Psi(\BF,\BA_{1},\BA_{2}) = \Psi_{\text{iso}}(\BF) + \Psi_{\text{ani}}(\BF,\BA_{1},\BA_{2})\,.
\end{equation}
Herein, the $\bar\BF$ is replaced with $\BF$ for the sake of convenience. \footnote{In the susequent treatment, the biological tissues will be assumed to behave perfectly incompressible.} Recall that, in the incompressible limit, the formulations (\ref{split1}) and (\ref{unsplit1}) lead to identical results. In this investigation, the volumetric part of the free energy function $U(J)$ will be excluded and the pressure term will be obtained from boundary conditions by enforcing exactly the incompressible deformation state under uniaxial and biaxial deformations. 
A canonical relation between the  stresses and the free energy function can be established in the sense
\begin{equation}
\BS = 2 \pa{\BC} \hat\Psi (\BC, \BA_1,\BA_2),~
\Btau = 2 \pa{\Bg} \tilde\Psi (\Bg, \BF, \BA_1,\BA_2) 
\label{eqn36}
\end{equation}
where $\BS$ and $\Btau$ are the second Piola-Kirchhoff and the Kirchhoff stresses, respectively. Thereafter, the Lagrangian moduli expression result from the second derivative of  $\hat\Psi$ with respect to $\BC$,
\begin{equation}
\IC = 4\pa{\BC\BC}^2 \hat\Psi (\BC, \BA_1,\BA_2)~.
\end{equation}
\subsubsection{Invariant-based anisotropic hyperelastic formulations}
The Lagrangian and Eulerian stresses for invariant-based formulations of anisotropic elasticity result from the application of the chain rule
\begin{equation}
\begin{array}{ll}
\BS = 2 \pa{\BC} \Psi =\Disp 2 \sum_{i=1} ^7 \pp{\Psi}{I_i}\pp{I_i}{\BC}\,\quad \text{and} \quad \Btau =  2 \pa{\Bg} \Psi=\Disp 2 \sum_{i=1} ^7 \pp{\Psi}{I_i}\pp{I_i}{\Bg}\,.
\end{array}
\label{eqntau}
\end{equation}
The closed form expressions for the Lagrangian and Eulerian stress tensors (\ref{eqntau}) require the following derivatives
\begin{equation}\label{derv-lag}
\begin{array}{lll}
\pa{\BC} I_1 &= \BG^{-1} \,,\quad
&\pa{\BC} I_2 = I_1 \BG^{-1} - \BC \,,\quad
\pa{\BC} J = \half J \BC^{-1}  \ ,\\[2Ex]
\pa{\BC} I_4 &= \BA_{1} \,,\quad
&\pa{\BC} I_5 = \BM_{1}\otimes\BC\BM_1 + \BM_1\BC\otimes\BM_{1} \ ,\\[2Ex]
\pa{\BC} I_6 &= \BA_{2} \,,\quad \quad
&\pa{\BC} I_7 = \BM_{2}\otimes\BC\BM_{2} + \BM_{2}\BC\otimes\BM_{2} \ .
\end{array}
\end{equation}
Substitution of (\ref{derv-lag}) into (\ref{eqntau})$_1$ gives the invariant-based representation for the second Piola-Kirchhoff stress tensor 
\begin{equation}
\begin{array}{rcl}
\BS &=& 2 (c_1 + I_1 c_2) \Bnone - 2 c_2  \BC  +  2 c_3 \BA_1 + 2 c_4 (\BM_1\otimes\BC\BM_1 + \BM_1\BC\otimes\BM_1)  \\[2Ex]
&+& 2 c_5 (\BM_2\otimes\BM_2) + 2 c_6 (\BM_2\otimes\BC\BM_2 + \BM_2\BC\otimes\BM_2)- p \BC^{-1} \,.
\end{array}
\label{stress-lag}
\end{equation}
The derivatives with respect to the spatial metric $\Bg$ follow from the push-forward
operation of (\ref{derv-lag})
\begin{equation}
\begin{array}{lll}
\pa{\Bg} I_1 &= \Bb \,, \quad
&\pa{\Bg} I_2 = I_1 \Bb - \Bb\Bg\Bb  \,, ~~
\pa{\Bg} J   = \half J \Bg^{-1}  \ ,\\[2Ex]
\pa{\Bg} I_4 &= (\Bm_1\otimes\Bm_1) \,, \quad \quad
&\pa{\Bg} I_5 = (\Bm_1\otimes\Bb\Bm_1+\Bm_1\Bb\otimes\Bm_1)  \ ,\\[2Ex]
\pa{\Bg} I_6 &= (\Bm_2\otimes\Bm_2) \,, \quad 
&\pa{\Bg} I_7 = (\Bm_2\otimes\Bb\Bm_2+\Bm_1\Bb\otimes\Bm_2) \ .
\end{array}
\label{derv-eul}
\end{equation}
Substitution of (\ref{derv-eul}) into (\ref{eqntau})$_2$  leads to the invariant-based representation for the second Kirchhoff stress expression
\eb
\begin{array}{rcl}
\Btau &=& 2 (\psi_1 + I_1 \psi_2)~ \Bb - 2 \psi_2 \Bb\Bg\Bb + 2 \psi_3 (\Bm_1\otimes\Bm_1) + 2 \psi_4 (\Bm_1\otimes\Bb\Bm_1+\Bm_1\Bb\otimes\Bm_1) \\[2Ex]
&+& 2 \psi_5 (\Bm_2\otimes\Bm_2) + 2 \psi_6 (\Bm_2\otimes\Bb\Bm_2+\Bm_2\Bb\otimes\Bm_2)- p \Bg^{-1}\,,
\end{array}
\label{stress-eul}
\ee
with $\Disp\psi_i =\pp{\Psi}{I_i}$ and $p=-J \pa{J}\Psi$. The stress expressions given in (\ref{stress-lag}) and (\ref{stress-eul}) can be expressed in terms of index notation for the principal directions
\eb
\begin{array}{r@{~}c@{~}l@{~}l@{~}l@{~}l@{~}l@{~}l@{~}l@{~}l@{~}l@{~}l@{~}l@{~}l}
S_i &=& 2 (\psi_1 + I_1 \psi_2) &-& 2 \psi_2  \lambda^2_i &+  2 \psi_3 M^2_{1i} &+& 4 \psi_4 M^2_{1i} \lambda^2_i &+& \Disp  2 \psi_5 M^2_{2i}  &+& 4 \psi_6 M^2_{2i} \lambda^2_i &-& \Disp\frac{1}{\lambda^2_i}p \,, \\[2Ex]
\tau_i &=& 2 (\psi_1 + I_1 \psi_2)\lambda^2_i &-& 2 \psi_2  \lambda^4_i  &+  2 \psi_3 m^2_{1i} &+& 4 \psi_4 m^2_{1i} \lambda^2_i&+&  2 \psi_5 m^2_{2i}   &+& 4 \psi_6 m^2_{2i} \lambda^2_i &-& p \,.
\end{array}
\label{ttti}
\ee
The Lagrangian moduli expression for invariant-based formulation can be written as
\begin{equation}
\begin{array}{c}
\IC = 4 \pa{\BC\BC} \Psi =\Disp 4 \sum_{i=1} ^7 \ppp{\Psi}{I_i}\ppp{I_i}{\BC}\,.
\end{array}
\label{eqntau}
\end{equation}
Furthermore, moduli expression can be additively decomposed into isotropic and anisotropic parts as follows
\begin{equation}
\IC=\IC_{iso}+\IC_{ani}~.
\end{equation}
The directional stiffness (DS) can be defined as
\begin{equation}\label{eqn:ds}
\textrm{DS} := \Bn\otimes\Bn:\IC_{ani}:\Bn\otimes\Bn\, \quad \text{with} \, \quad \Bn = cos\alpha \Be_1 + sin\alpha \Be_2~,
\end{equation}
where $\alpha \in[0,2\pi]$.
\subsubsection{Dispersion-type anisotropic formulations}
The tissue is considered a fiber-reinforced composite with the fibers distributed within an isotropic matrix in the dispersion-type anisotropic formulation; the models developed within this framework accurately describe the effect of the structural arrangement of the fibers on the mechanical response. Dispersion-type anisotropic  approaches utilize density distribution functions to represent the distributed fiber architecture of tissues. Let unit fiber direction $\Br$ on a unit-sphere be given in the undeformed configuration. The fiber density in direction $\Br$  is expressed with $\rho(\Br)$.  The unit fiber orientation vector can be expressed in spherical coordinates as follows
\begin{equation}
\Br =  \text{sin}\theta \text{cos}\phi \Be_1 +  \text{sin}\theta  \text{sin}\phi \Be_2 + \text{cos}\theta \Be_3\,,
\end{equation}
in terms of Eulerian angles $\theta\in[0,\pi]$ and $\phi \in[0,2\pi]$, see Figure \ref{sphere}. The Eulerian counterpart of $\Br$  is derived as $\Bt=\BF \Br$.
In this part, we will summarize the  two kinematic approaches for the modeling of dispersion-type anisotropy in soft biological tissues. \\[1ex] 
\inputfig{FIGs/sphere}{sphere}
\textbf{(i) Generalized structure tensor (GST) formulations:~} The generalized structure tensor is  defined as 
\begin{equation}
  \BnH=\frac{1}{\mid\SS\mid}\int_{\SS} \rho(\Br)\Br\otimes\Br \dA \quad
  \text{with} \quad \tr \BnH=1\,,
\label{GSTm}
\end{equation} 
and $\mid\SS\mid=4\pi$ for a unit sphere. Let $\Bh= \BF\BH\BF^T$ be defined as the Eulerian counterpart of the generalized structure tensor. The Lagrangian and Eulerian stresses can be obtained by applying the chain rule
\begin{equation}
\begin{array}{l@{\ =\ }c@{\ =\ }l}
\BS & 2 \pa{\BC} \psi & \Disp
2 \left[\pp{U}{J} \pp{J}{\BC} + \pp{\Psi_{\text{iso}}}{I_1} \pp{I_1}{\BC} +
\pp{\Psi_{\text{ani}}}{E} \pp{E}{\BC} \right]  ,
\end{array}
\label{eqngspk}
\end{equation}
\begin{equation}
\begin{array}{l@{\ =\ }c@{\ =\ }l}
\Btau & 2 \pa{\Bg} \psi & \Disp
2 \left[\pp{U}{J} \pp{J}{\Bg} +
\pp{\psi_{\text{iso}}}{I_1} \pp{I_1}{\Bg} +
\pp{\Psi_{\text{ani}}}{E} \pp{E}{\Bg} \right]  ,
\end{array}
\label{eqngtau}
\end{equation}
where
\begin{equation}
E := \BH : \BC -1 \qquad \text{or} \qquad E:=\Bh : \Bg -1 
\end{equation}
is the one-dimensional average anisotropic fiber strain akin to the Green-Lagrangian strain. The closed form expressions for the Lagrangian and Eulerian stress tensors require the following derivatives
\begin{equation}
\begin{array}{ccc}
\pa{\BC} E = \BH\,\qquad \text{and} \qquad \pa{\Bg} E = \Bh \,.
\end{array}
\label{Edc}
\end{equation}
Substituting (\ref{derv-lag}), (\ref{derv-eul}) and (\ref{Edc}) into (\ref{eqngspk}) and (\ref{eqngtau}), and by rearranging the terms, finally gives the generalized structure tensor-based representation of stresses for an anisotropic hyperelastic solid
\begin{equation}
\begin{array}{ll}
\BS = 2 \psi_1 \BG^{-1} + 2\psi_f \BH - p \BC^{-1} \,\quad\text{and} \quad \Btau = 2 \psi_1 \Bb + 2\psi_f \Bh - p \Bg^{-1} \,
\end{array}
\label{tg}
\end{equation}
with $\psi_f=\pa{E}{\Psi_{\text{ani}}}$. In the principal directions,~(\ref{tg}) can be written in terms of the principal stretches,
\begin{equation}
 \Disp S_i = 2 \psi_1 + 2 \psi_f H_i - \frac{1}{\lambda^2_i} p \,\quad \text{and} \quad \tau_i = 2 \psi_1 \lambda^2_i + 2 \psi_f h_i - p \,.
\end{equation}
The anisotropic part of the Lagrangian moduli expression can be written as
\begin{equation}
\IC_{ani} = 4 \ppp {\psi} {\BC} = 4 \psi_{ff} \BH_1\otimes\BH_1 + 4\psi_{ff} \BH_2\otimes\BH_2 \,\quad \text{with} \quad \psi_{ff}=\ppp{\Psi_{\text{ani}}} {\BE}~.
\end{equation}
\textbf{(ii) Angular integration (AI) formulations:~} The total free energy of fibers is calculated by 
\begin{equation}
\Psi_{\text{ani}}=\frac{1}{\mid\SS\mid}\int_{\mid\SS\mid} \rho(\Br)\psi_{fib}(\BF,\Br)\dA
\label{energyAI}
\end{equation}
where  $\SS$ represents a unit sphere with $\mid\SS\mid=4\pi$. 
Starting with the free energy (\ref{energyAI}), the Eulerian stress tensor follows by the application of the chain rule
\begin{equation}
\Btau = 2 \left[ \pp{\Psi_{\text{iso}}}{J}\pp{J}{\Bg} + \pp{\Psi_{\text{iso}}}{I_1}\pp{I_1}{\Bg} + \pp{\Psi_{\text{iso}}}{I_2}\pp{I_2}{\Bg} +  \pp{\Psi_{\text{ani}}}{\lambda_f} \pp{\lambda_f}{\Bg} \right]
\label{AIchain}
\end{equation}
Insertion of  (\ref{derv-eul}), $2 \pa{\Bg} \lambda_f = \lambda_f^{-1} \Bt\otimes\Bt$, and $\psi_f=\psi'_{fib}$ into (\ref{AIchain}) gives the Kirchhoff stress expression
\begin{equation}
\Btau = - p\Bg^{-1} + 2 (\psi_1 + I_1 \psi_2) \Bb -2 \psi_2 \Bb\Bg\Bb +  \int_{\SS} \rho(\Br) \lambda_f^{-1} \psi_f  \Bt\otimes\Bt \dA 
\label{tauAI}
\end{equation}
with $\psi_f=\pa{\lambda_F}{\psi_{fib}}$. The general form of the Kirchhoff stress based on the angular integration formulation is given in (\ref{tauAI}). In the principal directions, (\ref{tauAI}) can be written in terms of the principal stretches,
\eb
\begin{array}{rcl}
\tau_1 &=&\ds - p + 2(\psi_1 + I_1 \psi_2)\lambda^2_1 -2 \psi_2 \lambda_1^4  + \lambda_1^2 \int_{\SS} \rho(\Br) \lambda_f^{-1} \psi_f \text{cos}^2 \phi~ \text{sin}^3 \theta \dA  \,, \\[3Ex]
%
\tau_2 &=&\ds  - p + 2(\psi_1 + I_1 \psi_2)\lambda^2_2 -2 \psi_2 \lambda_2^4  + \lambda_2^2 \int_{\SS} \rho(\Br) \lambda_f^{-1} \psi_f \text{sin}^2 \phi~ \text{sin}^3 \theta \dA  \,, \\[3Ex]
%
\tau_3 &=&\ds - p + 2(\psi_1 + I_1 \psi_2)\lambda^2_3 -2 \psi_2 \lambda_3^4  + \lambda_3^2 \int_{\SS} \rho(\Br) \lambda_f^{-1} \psi_f \text{cos}^2 \theta~ \text{sin} \phi	\dA \,.
\label{tauAI1}
\end{array}
\ee
\inputfig{FIGs/DeformationModes}{deformartion_modes}
\subsection{Stresses under homogeneous deformations}
In this section, we briefly outline the deformation state of anisotropic biological tissues subjected to uniaxial and equibiaxial stresses. The deformation gradient tensor of anisotropic materials cannot be expressed in a straightforward manner as can be done for isotropic materials since the amount of contraction in the transverse direction and the shear deformations depend on the anisotropic fiber structure. Therefore, the deformation gradient tensor should be obtained iteratively to satisfy equilibrium conditions, e.g.~in case of uniaxial tension $\sigma_{22}=\sigma_{33}=\sigma_{12}=\sigma_{13}=\sigma_{23}=0$. To this end, we assume that the material is deformed in the principal directions corresponding to the principal axes of fiber orientations for AI- or GST-based models. For invariant formulations, the deformation axes are assumed to coincide with the symmetry axes for the two fiber family formulations, while in single fiber family formulations one of the principal deformation axes coincides with the fiber direction. Such an ansatz avoids \textsl{a priori}  shear strains/stresses, leading to diagonal deformation and stress tensors. In this context, the deformation gradient and the nominal stress expression under homogeneous uniaxial, equibiaxial or pure shear deformations can be expressed as follows
\begin{equation}
 \BF=\begin{bmatrix} F_{11} & 0 &0\\ 0 & F_{22} & 0 \\ 0 & 0 & F_{33}\end{bmatrix} 
  \quad \text{and} \quad
\BP=\begin{bmatrix} P_{1} & 0 &~0~\\ 0 & P_{2} & 0 \\ 0 & 0 & P_{3}\end{bmatrix}\,.
\end{equation}
\textbf{Uniaxial tension:} For an incompressible anisotropic hyperelastic solid, the deformation and stress states under uniaxial tension are
\begin{equation}
  \BF=\begin{bmatrix} \lambda & 0 &0\\ 0 & \frac{1}{\sqrt{\lambda}} & 0 \\ 0 & 0 & \frac{1}{\sqrt{\lambda}}\end{bmatrix} 
  \quad \text{and} \quad
\BP=\begin{bmatrix} P_{1} & 0 &~0~\\ 0 & 0 & 0 \\ 0 & 0 & 0\end{bmatrix} \,,
\end{equation}
see Figure \ref{def-modes}(a),(b). 
The first two invariants under uniaxial deformation read
\begin{equation}
I_1=\lambda^2+ \frac{2}{\lambda} \, \quad \text{and}\,\quad I_2=2\lambda+\frac{1}{\lambda^2} \,.
\end{equation}
For symmetrically orthotropic fibers, the equalities $I_4=I_6$ and $I_5=I_7$ hold leading to
\begin{equation}
I_4 = I_6 = \lambda^2 \text{cos}^2\varphi + \frac{1}{\lambda} \text{sin}^2\varphi \,\quad \text{and}\,\quad
I_5 = I_7 = \lambda^4 \text{cos}^2\varphi + \frac{1}{\lambda^2} \text{sin}^2\varphi \,,
\end{equation}
where $\varphi$ is the angle between the fiber and the symmetry axis $\Be_1$.
The components of the nominal stresses under uniaxial loading in the symmetry axes for the invariant-, generalized structure tensor- and angular integration-based formulations read

\begin{equation}\scalebox{0.97}{
\boxed{\begin{array}{r@{}r@{}l}
\text{Inv}: P_1&=&2 (\psi_1 + I_1 \psi_2)\lambda - 2 \psi_2  \lambda^3  +  2 \psi_3 m^2_{1\hskip 0.2ex 1}\lambda + 4 \psi_4 m^2_{1\hskip 0.2ex 1}\lambda^3 + \Disp  2 \psi_5 m^2_{2\hskip 0.2ex 1}\lambda  + 4 \psi_6 m^2_{2\hskip 0.2ex 1} \lambda^3- \frac{1}{\lambda}p \,,\\[2Ex]
\text{GST}: P_1&=& 2 \psi_1 \lambda+ 2 \psi'_f H_1 \lambda- \Disp \frac{1}{\lambda} p \,,\\[2Ex]
\text{AI}: P_1&=& n\lambda \int_{\SSs} \rho(\Br) \lambda_f^{-1} \psi_f \text{cos}^2 \phi \text{sin}^3 \theta	\dA + \Disp  2 \psi_1 \lambda  - \frac{1}{\lambda} p \,,
\end{array}}}
\end{equation}
\textbf{Equibiaxial tension:} For an incompressible hyperelastic anisotropic solid, the deformation and stress states under equibiaxial tension are
\begin{equation}
 \BF=\begin{bmatrix} \lambda & 0 &0\\ 0 & \lambda & 0 \\ 0 & 0 & \lambda^{-2}\end{bmatrix} 
  \quad \text{and} \quad
\BP=\begin{bmatrix} P_{1} & 0 &~0~\\ 0 & P_{2} & 0 \\ 0 & 0 & 0\end{bmatrix}\,.
\end{equation}
see Figure \ref{def-modes}(c). The first two invariants under equibiaxial deformation read
\begin{equation}
I_1=2\lambda^2+\frac{1}{\lambda^4}\,, \quad I_2=\lambda^4 + \frac{2}{\lambda^2}\,.
\end{equation}
The invariants associated with the structural tensors  are
\begin{equation}
I_4 = I_6 = \lambda^2 \text{cos}^2\varphi + \lambda^2 \text{sin}^2\varphi \,\quad \text{and}\,\quad
I_5 = I_7 = \lambda^4 \text{cos}^2\varphi + \lambda^4 \text{sin}^2\varphi \,.
\end{equation}
The components of the nominal stresses under equibiaxial loading in the symmetry axes for invariant-, generalized structure tensor- and angular integration-based formulations read
\begin{equation}\scalebox{0.96}{
\boxed{\begin{array}{rrl}
\text{Inv}: P_1&=&2 (\psi_1 + I_1 \psi_2)\lambda - 2 \psi_2  \lambda^3  +  2 \psi_3 m^2_{1\hskip 0.2ex 1}\lambda + 4 \psi_4 m^2_{1\hskip 0.2ex 1}\lambda^3 + 2 \psi_5 m^2_{2\hskip 0.2ex 1}\lambda  + 4 \psi_6 m^2_{2\hskip 0.2ex 1} \lambda^3- \Disp \frac{1}{\lambda}p \\[2Ex]
 P_2&=&2 (\psi_1 + I_1 \psi_2)\lambda - 2 \psi_2  \lambda^3  +  2 \psi_3 m^2_{2\hskip 0.2ex 1}\lambda + 4 \psi_4 m^2_{2\hskip 0.2ex 1}\lambda^3 + 2 \psi_5 m^2_{2\hskip 0.2ex 2}\lambda  + 4 \psi_6 m^2_{2\hskip 0.2ex 2} \lambda^3- \Disp \frac{1}{\lambda}p \\[2Ex]
\text{GST}: P_1&=& 2 \psi_1 \lambda+ 2 \psi'_f H_1 \lambda- \Disp \frac{1}{\lambda} p \\[2Ex]
 P_2&=& 2 \psi_1 \lambda+ 2 \psi'_f H_2 \lambda- \Disp \frac{1}{\lambda} p\\[2Ex]
\text{AI}: P_1&=& n\lambda \int_{\SSs} \rho(\Br) \lambda_f^{-1} \psi_f \text{cos}^2 \phi \text{sin}^3 \theta	\text{d}\theta\text{d}\phi+ 2 \psi_1 \lambda  - \Disp \frac{1}{\lambda} p \\[2Ex]
P_2&=& n\lambda \int_{\SSs} \rho(\Br) \lambda_f^{-1} \psi_f \text{sin}^2 \phi \text{sin}^3 \theta	\dA + 2 \psi_1 \lambda - \Disp \frac{1}{\lambda} p\,.
\end{array}}}
\end{equation}

\section{Hyperelastic material models}
\label{chp:b3}
In this section, we will review a common isotropic model and nine anisotropic hyperelastic models under two main categories: (i) strain invariant-based models, and (ii) fiber dispersion-based models. The free energy functions and necessary derivatives for the stress expressions are outlined. In order to emphasize the degree of anisotropy and the anisotropy distribution obtained from paramater identification process from each tissue is presented in terms of polar plots for each constitutive model investigated. The polar plots consider the density distribution (DD) function for the plane of interest along with the directional stiffness (DS) computed from equation (\ref{eqn:ds}). For the invariant based formulations, the polar plots of the directional stiffness obtained merely from the anisotropic part of the initial elastic moduli tensor is used for visualization of the plane of interest. 
\subsection{Invariant-based models}
\subsubsection{The neo-Hookean model}
The invariant-based models assume perfect alignment of fibers embedded into an isotropic ground matrix. Models have been presented for two families of fibers symmetrically disposed.
The neo-Hookean model is the most fundamental model for hyperelastic constitutive models. Many researchers represent the isotropic ground matrix of soft tissues with the neo-Hookean model~\cite{weiss1996,holzapfel2000,Holzapfel2005,Alastrue2009,Alastrue2010}. Therefore, it deserves a separate description. Based on Wall's treatment of elasticity of a molecular network, Treloar~\cite{Treloar1943} proposed the following free energy function
\begin{equation}
  \Psi_{\text{iso}}=\frac{1}{2}n k_B \theta\left(\lambda_{1}^{2}+\lambda_{2}^{2}+\lambda_{3}^{2}-3\right)\quad \text{or}\quad
\Psi_{\text{iso}}=\frac{1}{2}\mu\left(I_1-3\right)
\label{neo1}
\end{equation}
with $\mu=nk_B\theta$. Herein, $n$ is the volume specific chain density, $k_B$ is the Boltzmann constant and $\theta$ is the absolute temperature, and $\mu=nk_B\theta$ is the shear modulus. The nonzero derivative of $\psi$ with respect to the invariants is
\begin{equation}
  \psi_1=\frac{\mu}{2}~. 
\end{equation}
\subsubsection{Newman-Yin (NY) model}
\inputfig{FIGs/polarplotforinv/ny}{ny_polar} 
Newman and Yin~\cite{Newman1998} proposed a exponential free energy form analogous to the one proposed by Fung et al.~\cite{Fung1979} to describe hyperelastic behavior of mitral valve tissue. They assumed that the material is transversely isotropic and that the free energy function depends on the first and the fourth invariants. They observed that for a constant $I_4$, both $\psi_1$ and $\psi_4$ increase nonlinearly. Therefore, Newman and Yin proposed, 
\begin{equation} \label{free-ny}
\Psi_{\text{ani}} = k_0 (\textrm{exp}(Q)-1) \quad \text{with} \quad Q = (k_1 (I_1 -3)^2 + k_2 (\sqrt{I_4} -1)^4)+ k_2 (\sqrt{I_6} -1)^4)~ \,.
\end{equation}
$Q$ is the quadratic function of the invariants. In order to have a strain energy increase with increasing $I_4$, $k_0$ should be positive. Also, if the tissue is not able to support compressive load, $k_2$ should be positive. The original model was proposed for tissues with a single family of fibers. In this work, we have extended the formulation considering two families of fibers as given in (\ref{free-ny})$_2$ by incorporating the latter term depending on $I_6$. The derivatives of $\Psi_{\text{ani}}$ with respect to the invariants are 
\begin{equation}
\begin{array}{lll}
  \psi_1 &=& 2 k_0 k_1 (I_1 -3)\textrm{exp}(Q) \,, \quad \psi_4 = 2 k_0 k_2 \textrm{exp}(Q) \Disp \frac{1}{\sqrt{I_4}} (\sqrt{I_4}-1)^3\,,\\[2Ex] 
\psi_6 &=& 2 k_0 k_2 \textrm{exp}(Q) \Disp \frac{1}{\sqrt{I_6}} (\sqrt{I_6}-1)^3 \,.
\end{array}
\end{equation}
The polar plots of the directional stiffness for the NY model obtained with the identified parameters for each tissue are depicted in Figure~\ref{ny-polar}.
\subsubsection{Holzapfel-Gasser-Ogden (HGO) model}
Holzapfel et al.~\cite{holzapfel2000} proposed a constitutive model for arteries. They additively decomposed the isochoric free energy function into isotropic and anisotropic parts,
\begin{equation}
\Psi = \Psi_{\text{iso}}(I_1) + \Psi_{\text{ani}}(I_4,I_6) \,.
\end{equation}
They utilized the neo-Hookean model (\ref{neo1}) for the isotropic part since collagen fibers are thought to not contribute to the mechanical behavior of the tissue at low pressures. The free energy stored by two families of collagen fibers is described as 
\begin{equation}
\Psi_{\text{ani}}(I_4,I_6)=  \frac{k_1}{2k_2}\sum_{i=4,6} \left(\text{exp}[k_2 \langle I_i -1\rangle^2]-1\right)
\end{equation}
where $k_1 > 0$ is a stress-like parameter and $k_2 > 0$ is a dimensionless parameter. The collagen fibers do not support compressive stresses due to their wavy nature. Therefore, the Macauley brackets 
\begin{equation}
\langle (\bullet) \rangle = \frac{(\bullet) + \mid(\bullet)\mid}{2}
\end{equation}
are utilized in order to filter out the tensile stretches. The derivatives of $\psi_i$ with respect to the $i^\text{th}$ invariants are 
\begin{equation}
\begin{array}{rll}
  \psi_1 &=& \Disp \frac{\mu}{2}\,,\quad  \psi_4 =  k_1 \langle I_4-1\rangle\textrm{exp}\left(k_2 \langle I_4 -1\rangle^2-1\right)\,,\\[2Ex]   \psi_6 &=& k_1 \langle I_6-1\rangle  \textrm{exp}\left(k_2 \langle I_6 -1 \rangle^2-1\right)~.
\end{array}
\end{equation}
The polar plots of the directional stiffness for the HGO model obtained with the identified parameters for each tissue are depicted in Figure~\ref{hgo-polar}.
\inputfig{FIGs/polarplotforinv/hgo}{hgo_polar} 
\subsubsection{Holzapfel-Sommer-Gasser-Regitnig (HSGR) model}
Holzapfel et al.~\cite{Holzapfel2005} proposed a free energy function of the form $\Psi(I_1,I_4) = \Psi_{\text{iso}}(I_1) + \Psi_{\text{ani}} (I_1,I_4)$. The model uses the  neo-Hookean model (\ref{neo1}) as the isotropic part of the free energy function. The anisotropic part of the free energy function has a mixed representation 
\begin{equation}
\Psi_{\text{ani}} =  \text{sgn}\langle I_4-1\rangle\frac{k_1}{k_2}\left(\textrm{exp}\{k_2[(1-p)( I_1 -3 )^2 + p\langle I_4-1 \rangle ^2]\}-1\right)~,
\label{enhsgr}
\end{equation}
\inputfig{FIGs/polarplotforinv/hsgr}{hsgr_polar} 
where $k_1 > 0 $ is a stress-like parameter and $k_2 >$ is a dimensionless parameter. The measure of anisotropy parameter $p \in [0,1]$ interpolates between the contributions of the first $I_1$ and fourth invariants $I_4$. It mimics the degree of fiber dispersion phenomenologically.  The switch function sgn$\langle I_4-1\rangle$ enforces the tension-only condition and is activated for $I_4-1>0$.  The model recovers the HGO model for $p=1$. The fiber related terms drop and the model reduces to an exponential isotropic model for $p=0$. The nonzero derivatives of the free energy function are
\begin{equation}
\begin{array}{rl}
  \psi_1 =& \ds\frac{\ds\mu}{2} + 2 k_1 (1-p) (I_1 -3)\text{sgn}\langle I_4-1\rangle \textrm{exp}\{k_2[(1-p)(I_1 -3)^2 + p\langle I_4-1 \rangle^2]\}\,,  \\[1.5Ex]
  \psi_4 =& 2 \ds\frac{k_1}{k_2} p \langle I_4 -1 \rangle \textrm{exp}\{k_2[(1-p)(I_1 -3)^2 + p\langle I_4-1 \rangle^2]\}~.
\end{array}
\end{equation}
\inputfig{FIGs/polarplotforinv/os}{os_polar} 
\subsubsection{Ogden-Saccomandi (OS) model}
Ogden and Saccomandi~\cite{Ogden2007}, following the work of  Horgan and Saccomandi~\cite{Horgan2005},  proposed a logarithmic constitutive law for arterial tissue with two fiber families in which the fiber extension is limited. They additively decomposed the free energy function into isotropic and anisotropic parts,
\begin{equation}
\Psi = \Psi_{\text{iso}}(I_1) + \Psi_{\text{ani}}(I_4,I_6) \,.
\end{equation}
To model the isotropic behavior of the tissue, they adapted the well-known rubber elasticity model of Gent \cite{Gent1996}. The free energy function of the Gent model reads 
\begin{equation}
\Psi_{\text{iso}} = -\frac{1}{2} \mu J_m \ln\left(1-\frac{I_1-3}{J_m}\right) \,,
\label{gent}
\end{equation}
where $\mu$ is the shear modulus and $J_m$ is the parameter that controls the chain extensibility limit for the matrix material. The deformation limit for the first invariant is $I_1 < 3 + J_m$ and the stresses tend to infinity asymptotically at this limit. As $J_m  \rightarrow \infty$, the isotropic part of the free energy function (\ref{gent}) recovers the neo-Hookean model. A similar model of Gent \cite{Gent1996} was proposed by Horgan and Saccomandi~\cite{Horgan2005} for transversely isotropic materials. Instead of limiting the polymer chain extensibility, the model of Horgan and Saccomandi limits the extensibility of the fibers.  The anisotropic part of the free energy function is
\inputfig{FIGs/dispersionparameterGasser}{measure-dispersion} 
\begin{equation}
\Psi_{\text{ani}} = - \frac{k_1}{2}J_f \sum_{\alpha = 4,6} \ln \left(1-\frac{\langle I_{\alpha}-1 \rangle^2}{J_f} \right) \,,
\label{ogdenani}
\end{equation}
where $k_1$ is a stress-like parameter, and $J_f$ is the limiting parameter of extensibility of collagen fibers. The constraint
\begin{equation}
I_{\alpha} < \sqrt{J_f} + 1 \,,\quad \alpha=4,6 \,.
\end{equation}
confines the stretches beyond the fiber extensibility limit. The derivatives of the free energy function are
\begin{equation}
\begin{array}{rll}
  \psi_1 &=& \Disp \frac{\mu}{2}\left(\frac{J_m}{J_m - I_1 +3}\right) \,,\quad  \psi_4 = k_1 J_f \left(\Disp \frac{\langle I_4 -1 \rangle}{J_f-\langle I_4-1\rangle^2} \right)\,,\\[2Ex]   
 \psi_6 &=& k_1 J_f \left( \Disp\frac{\langle I_6 -1 \rangle}{J_f-\langle I_6-1\rangle^2}\right)~.
\end{array}
\end{equation}
The polar plots of the directional stiffness for the OS model obtained with the identified parameters for each tissue are depicted in Figure~\ref{os-polar}.
\subsection{Dispersion-type anisotropic constitutive models}
It has been determined that not only do the mean orientation of fibers within soft tissues affect their mechanical properties, but also the amount by which they are dispersed around this mean~\cite{holzapfel2019}. For instance, if the majority of fibers are orientated along the direction of the mean, the behavior will be very different to that if the fibers are distributed in a cone-like span around the mean orientation~\cite{Gasser2006}. For this reason, a number of models have been developed to incorporate this dispersion into the constitutive law describing the stress-strain relation of the tissue. These fiber dispersion-based models utilize a probability distribution function to model the histological structure of tissues. The models which have an angular integration approach and a generalized structure tensor approach are outlined in this section.
\subsubsection{GST-based Gasser-Ogden-Holzapfel (GOH) model}
Gasser et al.~\cite{Gasser2006}  assumed that the tissue's fibers are distributed rotationally symmetric around a mean fiber orientation direction $\BM$. They utilized a planar $\pi$-periodic von-Mises distribution as a fiber density distribution around $\BM$. The von-Mises distribution function is a one-dimensional probability distribution which is a function of $\Theta$ and concentration parameter $b$. The standard von-Mises distribution function is
\begin{equation}
\bar\rho(\Theta)=\frac{\exp[b(\cos(2\Theta))}{2\pi I_0(b)} \quad \text{with} \quad
I_0(b)=\frac{1}{\pi}\int^{\pi}_{0} \exp(b\cos\Theta)\text{d}\Theta\,,
\label{standardvonM}
\end{equation}
where $I_0(b)$ is a modified Bessel function of the first kind of order zero. Applying the normalization condition to (\ref{standardvonM}) gives the relation
\begin{equation}
I \equiv \int_0 ^{\pi} \bar\rho(\Theta)\text{sin}\Theta\text{d}\Theta \equiv \frac{\text{exp}(-b)}{2\sqrt{2\pi b}}\frac{\text{erfi}(\sqrt{2b}}{I_0(b)} \,.
\end{equation}
The normalized von-Mises distribution is
\begin{equation}
\rho(\Theta) = \frac{\bar\rho(\Theta)}{I} = 4\sqrt{\frac{b}{4\pi}}\frac{\text{exp}\left(b\text{cos}(2\Theta)+1\right)}{\text{erfi}\sqrt{2b}}  \,.
\label{eq-vonmises}
\end{equation}
By inserting the von-Mises type density distribution function (\ref{eq-vonmises}) into (\ref{GSTm}), the generalized structure tensor can be written as 
\begin{equation}
\BH=\kappa\Bnone+(1-3\kappa)\BM\otimes\BM \,\quad \text{where}\,\quad\kappa=\frac{1}{4}\int^{\pi}_{0} \rho(\Theta)\sin^{3}\Theta \text{d}\Theta
\end{equation}
\inputfig{FIGs/polarplots/goh_polar}{goh_polar} 
is the fiber dispersion parameter. A one-to-one relation exists between the dispersion parameter $\kappa$ and the concentration parameter $b$, see Figure \ref{measure-dispersion}. Hence, $\kappa\in[0,1/3]$ enters the constitutive model as an additional material parameter responsible for the degree of dispersion. The lower limit $\kappa=0$ recovers the invariant-based anisotropy and the upper limit $\kappa=1/3$ leads to an isotropic constitutive response. The polar plots of the density distribution and directional stiffness for the GOH model obtained with the identified parameters for each tissue are depicted in Figure~\ref{goh-polar}. Gasser et al. ~\cite{Gasser2006} can be considered as the GST counterpart of the GHO model ~\cite{holzapfel2000}. The anisotropic part of the free energy function reads
\begin{equation}
\Psi_{\text{ani}}(\BC,\BH_i)=\frac{k_1}{2k_2}[\text{exp}(k_2 E^2_i)-1]\,,\quad i=1,2
\label{eq3.6}
\end{equation}
where $E_i=\BH_i:\BC- 1$ replaces the fourth and sixth invariants in the GHO model. The classical neo-Hookean model (\ref{neo1}) is utilized for the isotropic part of the free energy function.
The derivatives of the free energy function are
\begin{equation}
\psi_1 = \mu\,\quad \text{and} \quad \psi_{f} =  k_1 E_i \text{exp} (k_2 E_i ^2)\,.
\label{eq3.7}
\end{equation}
\inputfig{FIGs/bivariate-dispersion}{bivariate} 
\subsubsection{GST-based Holzapfel-Niestrawska-Ogden-Reinisch-Schriefl (HNORS) model}
Holzapfel et al.~\cite{Holzapfel2015b} take into account both the in- and out-of-plane dispersion of fibers based on the observations of Schriefl et al.~\cite{Schriefl2012a,Schriefl2012b,Schriefl2012,Schriefl2013} where they recorded that the fibers are dispersed both in-plane and out-of-plane. Their observations reveal that no correlation exists between in- and out-of-plane dispersions. Based on these arguments, the probability density function is multiplicatively decomposed as
\begin{equation}
\rho(\Br)=\rho_{ip}(\Phi)\rho_{op}(\Theta).
\end{equation}
For in-plane distribution, they considered a basic von-Mises distribution
\begin{equation}
  \rho_{ip}(\Phi)=\frac{\exp[a(\cos(2\Theta))}{I_0(a)}
    \quad\text{with}\quad
    I_0(a)=\frac{1}{\pi}\int^{\pi}_{0} \exp(x\cos\alpha)d\alpha
\end{equation}
where $a$ is the concentration parameter and $ I_0(a)$ is the modified Bessel function of the first kind of order zero. The out-of-plane distribution is in the form 
\begin{equation}
\rho_{op}(\Theta)=2\sqrt{\frac{2b}{\pi}}\frac{\exp[b(\cos(2\Theta))-1]}{\text{erf}(\sqrt{2b})}~.
\end{equation}
The measures of dispersion in the in-plane and the out-of-plane directions read
\begin{equation}
\kappa_{ip}=\frac{1}{\pi}\int^{\pi}_{0} \rho_{ip}(\Phi)\sin^{2}\Phi \text{d}\Phi\,\quad \text{and} \quad \kappa_{op}=\int^{\pi/2}_{0} \rho(\Theta)\sin^{3}\Theta \text{d}\Theta
\end{equation}
The structure tensor $\BH$  has the form
\begin{equation}
\BH=2\kappa_{ip}\kappa_{op}\Bnone + 2\kappa_{op}(1-2\kappa_{ip})\BM_i\otimes\BM_i + (1-2\kappa_{op}-2\kappa_{ip}\kappa_{op})\BM_n\otimes\BM_n
\end{equation}
where $\BM_i$ is the in-plane mean fiber direction, whereas $\textbf{M}_n$ is the out-of-plane vector. The Eulerian counterpart of the generalized structure tensor reads
\begin{equation}
\Bh=2\kappa_{ip}\kappa_{op}\Bnone + 2\kappa_{op}(1-2\kappa_{ip})\Bm_i\otimes\Bm_i + (1-2\kappa_{op}-2\kappa_{ip}\kappa_{op})\Bm_n\otimes\Bm_n 
\end{equation}
where $\Bm_i=\BF\BM_i$ and $\Bm_n=\BF\BM_n$, respectively. The lower and upper bounds for the in-plane and out-of-plane dispersion parameters are, respectively,
\begin{equation}
 0 \le\kappa_{op}\le 1/2 \qquad \textrm{and} \qquad 0 \le\kappa_{ip}\le 1~.
\end{equation}
\inputfig{FIGs/polarplots/hnors_polar}{hnors_polar} 
The polar plots of the density distribution and directional stiffness for the HNORS model obtained with the identified parameters for each tissue are depicted in Figure~\ref{hnors-polar}.
The free energy functions are identical to the GOH model\cite{Gasser2006} as depicted
in (\ref{eq3.6}).
\subsubsection{AI-based Alastru{\'{e}} -Martinez-Doblar{\'{e}} -Menzel (AMDM) model}
The angular integration-based anisotropic model of Alastru{\'{e}}  et al.~\cite{Alastrue2009} takes into account rotationally symmetric fiber dispersion based on the micro-sphere model. The model utilizes a planar $\pi$-periodic von-Mises distribution for the fiber density distribution $\rho(\Br,\BM)$ around a mean direction $\BM$, in the same sense as Gasser et al.~\cite{Gasser2006}. $\Br$ is the unit orientation vector of a micro-fiber and $\Bt:=\BF\Br$ is the Eulerian counterpart of the Lagrangian fiber vector. The affine-stretch of a single fiber in the orientation direction $\Br$ reads
\begin{equation}\label{fiber-stretch}
\lambda_f:=\sqrt{\Bt_\flat\cdot\Bt}\, \quad \text{where} \quad\Bt_\flat:=\Bg\Bt~.
\end{equation}
The macroscopic free energy corresponding to one family of fibers with the mean direction $\BM$ and with n fibers per unit volume is defined as 
\begin{equation}\label{cont-ai}
\Psi_{\text{ani}}(\Bg,\BF) =\big\langle \rho\psi_{fib}(\lambda_f) \big\rangle=\frac{1}{\mid\SS\mid} \int_{\SS} \rho(\Br;\BM)\psi_{fib}(\lambda_f)dA
\end{equation}
where $\psi_{fib}$ is the free energy function associated with the orientation direction $\Br$ and $\mid\SS\mid=4\pi$ for a unit-sphere. 
For the isotropic ground matrix, they utilized the neo-Hookean free energy function. The anisotropic free energy function is as follows
\inputfig{FIGs/polarplots/amdm_polar}{amdm_polar} 
\begin{equation}
  \psi_f=\begin{cases}
    0\,, & \text{for $\lambda_f<1$}\\\ds
    \frac{k_1}{2k_2}[\exp(k_2[\lambda_f^2-1]^2)-1]\,, & \text{for $\lambda_f\geq1$}~.
  \end{cases}
\end{equation}
The contribution of each family of fibers to the macroscopic isochoric Kirchhoff stresses can be expressed as a continuous average, including the orientation distribution function, namely
\begin{equation}\label{tau-AMDM}
\begin{array}{cc}
\tau_f = \big\langle \rho\psi_f\bar{\lambda}^{-1}\Bt\otimes\Bt\big\rangle
\quad
\text{where} 
\quad
\psi_f:= \ds\pp{\psi_{fib}}{\lambda_f} = 2k_1\lambda_f[\lambda^2_f -1]\exp(k_2[\lambda^2_f -1]^2) \,.
\end{array}
\end{equation}
The continuous average in (\ref{cont-ai}) and (\ref{tau-AMDM}) is approximated by
\begin{equation}
\langle (\bullet) \rangle =\frac{1}{\mid\SS\mid}\int_{\Omega} (\bullet)dA \approx \sum_{i=1}^{m} w^i \left(\bullet\right)^i
\end{equation}
where ${ w^i }_{i=1,...,m}$ are the weight factors associated with the discrete orientation directions ${\Br^i }_{i=1,...,m}$. The polar plots of the density distribution and the directional stiffness for the AMDM model obtained with the identified parameters for each tissue are depicted in Figure~\ref{amdm-polar}. The rotationally symmetric $\pi$-periodic, normalized von-Mises distribution reads 
\begin{equation}
\rho(\Theta) = \frac{\bar\rho(\Theta)}{I} = 4\sqrt{\frac{b}{4\pi}}\frac{\text{exp}\left(b\text{cos}(2\Theta)+1\right)}{\text{erfi}\sqrt{2b}}  \,.
\label{eq-vonmisesAI}
\end{equation}
In this regard, the AMDM model can be considered as the AI counterpart of the GOH model.
\subsubsection{AI-based Alastru{\'{e}} -Sa{\'{e}}z-Martinez-Doblar{\'{e}}  (ASMD) model}
As an extension of their previous model, Alastru{\'{e}}  et al.~\cite{Alastrue2010} included the Bingham distribution in their constitutive model. This distribution function exhibits andipodal symmetry and is expressed as
\begin{equation}\label{bingham-dd}
  \rho(\Br;\BZ,\BQ)
  =[F_{0000}(\BZ)]^{-1}\text{etr}(\BZ\cdot\BQ^T\Br\cdot\Br^T\BQ)
\end{equation}
where $\text{etr}(\bullet)\equiv\exp(\tr(\bullet))$, $\BZ$ is a diagonal matrix with eigenvalues $[\kappa_1,\kappa_2,\kappa_3]$, $\textbf{Q}$ is orthogonal orientation matrix such that $\BA=\BQ\cdot\BZ\cdot\BQ^T$ and $F_{0000}(\BZ)$ is defined as 
\begin{equation}
F_{0000}(\BZ)=[4\pi]^{-1}\int_{\SS}(\text{etr}(\BZ:\Br\Br^T)\dA = {}_1F_1(\frac{1}{2};\frac{2}{3};\BZ)
\end{equation}
\inputfig{FIGs/polarplots/asmd_polar}{asmd_polar} 
where ${}_1F_1$ is a confluent hypergeometric function of the matrix argument. The shape of the distribution is controlled by $\kappa_1$, $\kappa_2$, and $\kappa_3$.  In this regard, the ASMD model is a slight modification of the AMDM model that utilizes the Bingham distribution, where the density distribution given in (\ref{bingham-dd}) replaces $\pi$-periodic von-Mises distribution  (\ref{eq-vonmisesAI})  in the equations (\ref{cont-ai}, \ref{tau-AMDM}). The polar plots of the density distribution and the directional stiffness for the ASMD model obtained with the identified parameters for each tissue are depicted in Figure~\ref{asmd-polar}.
\inputfig{FIGs/polarplots/dbb_polar}{dbb_polar} 
\subsubsection{AI-based Driessen-Bouten-Baaijens (DBB) model}
Driessen et al.~\cite{Driessen2005} presented an extended version of the Holzapfel et al.~\cite{holzapfel2000} (HGO) model which included a fiber volume fraction. They applied the rule of mixtures and expressed the isochoric Kirchhoff stress for multiple fiber directions as follows
\begin{equation}
\Btau=\Btau^m+\frac{1}{4\pi}\int_{\SS} \frac{v_f}{\lambda_f^2}(\tau_f-\frac{1}{\lambda_f^2}\Bt\cdot\Btau^m\Bt)\Bt\otimes\Bt \dA
\end{equation}
where $\Btau^m$ is the isotropic matrix stress, $v_f$ is the volume fraction of fibers and $\tau_f$ is the fiber stress for a given orientation direction $\Br$. The isotropic matrix material is modeled as a neo-Hookean material with a shear modulus $\mu$. The stress expressions for the isotropic matrix and fibers are
\begin{equation}
\Btau^m=\mu(\Bb- \Bnone)\,\quad \text{and} \quad \tau_f = k_1 \lambda_f^2\left[k_2\text{exp}\left(\lambda_f^2 -1\right)-1\right]~,
\end{equation}
where the fiber stretch $\lambda_f$ is given in (\ref{fiber-stretch}). The DBB model utilizes a planar Gaussian distribution around the mean fiber orientation for the orientation fiber volume content
\begin{equation}
v_f(\phi)=A \bar v_f(\phi) \quad \text{with} \quad \bar v_f(\phi):=\exp[\frac{-(\phi-\vartheta)^2}{2\sigma^2}]
\end{equation}
where $\vartheta$ is the mean value, $\sigma$ is the standard deviation and A is the normalization constant which is defined as
\begin{equation}
\text{A}=\frac{v_{tot}}{\ds\frac{1}{\mid\SS\mid}\int_{\SS}\Disp \bar v_f(\phi)\dA} \,.
\end{equation}
Therein, $v_{tot}$ is the total fiber volume fraction.  Unlike the AMDM and ASMD models, the DBB model excludes the volume fraction of the matrix for each orientation direction leading to a structure
\begin{equation}
\frac{1}{\lambda_f^2}\Btau:\Bt\otimes\Bt=(1-v_f)\frac{1}{\lambda_f^2}\Btau^m:\Bt\otimes\Bt +  v_f\tau_f
\end{equation}
where, the directional stress additively decomposes into the matrix contribution and collagen fibers contribution proportional to their volume fraction in each orientation direction, respectively. The polar plots of the density distribution for the DBB model obtained with the identiﬁed parameters for each tissue are depicted in Figure~\ref{dbb-polar}.

\section{Experimental Studies}
\label{chp:b4}
In this section, we outline the experiments that are used for the comparative investigations made in this paper. To this end, three distinct human tissues;  (i) the aneurysmatic abdominal aorta, (ii) the linea alba, and (iii) the anterior rectus sheath are selected for the subsequent analysis. Histologically, the AAA tissue has two alternating families of collagen fibers each orientated at an angle approximately $30^{\circ}$ to the circumferential direction \cite{niestrawska2016}, the linea alba has collagen fibers situated in the transverse direction \cite{cooney2016}, while the rectus sheath's principal fiber direction is in the longitudinal direction \cite{Martins2012}. In the experimental studies, the aneursymatic arterial tissue is subjected to equibiaxial tension, whereas the data of the linea alba and anterior rectus sheath originate from uniaxial tension tests executed in two mutually orthogonal directions. In this part, the experimental protocol is briefly reviewed.

\subsection{Uniaxial tension}
\noindent Uniaxial tensile tests are performed with samples having a length-to-width ratio of at least 4:1, see Figure \ref{def-modes}(a). To experimentally investigate anisotropy, samples are usually cut and tested along both the $\Be_1$- and $\Be_2$-directions, as seen in Figure \ref{def-modes}, with the assumption that no anisotropy derives from the $\Be_3$-direction. As the reference axis can change depending on the samples being tested, all following discussion will take place in the context of the $\Be_1$-direction samples, with the understanding that the same principles can be easily adapted by the reader for the $\Be_2$-direction samples. 
\noindent Typically for the uniaxial tension tests of metals or polymers, dogbone samples are employed to reduce the influence of end effects, and encourage fracture within the middle of the specimen and not at the grip location. However, for soft biological tissues, rectangular samples are often used instead due to the soft and delicate nature of the material, rendering the punching of precise dogbone samples generally very difficult.  \\

\noindent The first study selected to present uniaxial tension data on soft biological tissues was by Cooney et al. \cite{cooney2016}, who investigated the uniaxial and equibiaxial tensile behavior of the human linea alba. Cooney et al. \cite{cooney2016} obtained 13 freshly frozen human cadaveric abdominal walls for their study and from these extracted the linea alba for their experiments. The linea alba is a collagenous part of the ventral abdominal wall whose collagen fibers are known to be anisotropically arranged. Prior to sample extraction, the human cadaveric abdominal walls were allowed to defrost for 36 hours at $4^{\circ}$C. The linea alba were then extracted from each abdominal wall and stored in a phosphate-buffered saline solution until testing. For the uniaxial tension tests, 14 rectangular samples in total were obtained from 7 of the linea alba specimens. Seven samples were cut in both the longitudinal and transverse directions, with longitudinal and transverse referring to their respective anatomical axes. The length-to-width ratio of the rectangular samples was as close to the uniaxial test condition as possible for the size of available tissue, with this equaling approximately 2:1 (length:width) for both directions. Prior to testing, samples were mounted in grips lined with emery paper to reduce any slippage during loading. The grips were tightened using a 0.2 N-m torque wrench to prevent over-tightening of the grips causing damage to the tissue. Six black dots were applied to the surface of the samples
, before securing them in the modular mechanical testing machine which was fitted with a 300 N load cell. During testing, a high definition camera was used to collect images at a rate of 2 Hz. The post-processing of these images allowed for the calculation of the sample strains by analyzing the deformation of each black dot throughout the stretching. A prestress of 0.1 MPa was applied to all samples to remove any slack present. The tests were conducted at a quasi-static strain rate of 28.5\% min$^{-1}$.
The transverse direction was found to be much stiffer than the longitudinal direction, displaying the anisotropy of the tissue. Cooney et al. \cite{cooney2016} found the slope of the most linear section of the curve to be approximately 72 MPa in the transverse direction and 8 MPa in the longitudinal direction. The second study selected was by Martins et al.\cite{Martins2012}, who investigated the uniaxial tensile behavior of the human anterior rectus sheath. Martins et al. extracted the tissue samples from 12 fresh female cadavers. Six rectangular samples were cut both parallel (longitudinal) and perpendicular (transverse) to the fiber direction with a length-to-width ratio of 4:1. The thickness of the samples were 1.00 $\pm$ 0.17 mm. A combination of sandpaper and Velcro tape were used between the grips to stop the soft tissue sample from slipping during testing. A 200 N load cell was fitted to the traction machine to obtain the force, and the displacement of the samples was measured using a displacement sensor. The uniaxial tests were carried out at room temperature and were conducted until rupture at a rate of 5 mm min$^{-1}$. The tensile behavior was observed to be nonlinear in both the fiber direction and transverse to the fiber direction. 
\subsection{Equibiaxial tension}
\noindent Equibiaxial tensile tests are utilized in the characterization of soft biological tissues as the test condition is often considered to be more representative of the stresses experienced {\it in vivo} \cite{cooney2016}. Square samples are usually preferred for the tests and the samples are equally loaded under tension on all four sides, as seen in Figure \ref{def-modes}(c). Due to this, the principal stretch in the $\Be_1$- and $\Be_2$-directions are equal to the experimental stretch, i.e. $\lambda_1 = \lambda_2 = \lambda$, leaving the $\Be_3$-direction the only one unconstrained. If the material is considered incompressible, and due to the assumption of symmetry, $\lambda_3 = \lambda^{-2}$. 
\noindent The equibiaxial tension study selected was an investigation by Niestrawska et al. \cite{niestrawska2016} on human abdominal aortas. They conducted equibiaxial tension experiments on both healthy and aneurysmatic abdominal aortas. The data used here was extracted from the experimental results on the abdominal aortic aneurysms which were obtained through surgical procedures conducted at the Department of Vascular Surgery, Medical University Graz, Austria. Aneurysms are defined as irreversible, localized dilatation of a vessel, which can in some cases lead to complete wall rupture. Eleven wall samples in total were collected from open aneurysm repair at the anterior side. The samples were stored at $4^{\circ}$C in Dublecco's modified Eagle's medium prior to testing.  To prepare for testing, 20 $\times$ 20mm patches were cut from the aneurysmatic abdominal aorta walls. Some were large enough for two test specimens to be prepared from a single extracted sample. The specimens were tested with their layers intact as the authors found a clear separation of layers of the aneurysms to be impossible. At this point, the thickness of the specimens were measured. A scatter pattern of black dots was then applied to the surface of the specimens using a spray to allow the displacements to be optically measured. To mount the specimens in the biaxial testing machine, hooked surgical sutures were used. The specimens were then submerged in a 0.9\% physiological saline bath which was heated up to 37 $\pm$ $0.1^{\circ}$C. The tests were conducted using a stretch-driven protocol, starting at 2.5\% deformation and increasing in 0.025 stretch steps until rupture. This was carried out at a quasi-static strain rate of 3mm min$^{-1}$, and with a stretch ratio of $\lambda_{axial}:\lambda_{circ} = 1:1$, where $\lambda_{axial}$ is the stretch in the axial direction and $\lambda_{circ}$ is the stretch in the circumferential direction. Four preconditioning cycles were conducted after each increase in step, with the fifth cycle being recorded for data analysis. 
\noindent The results presented in the study by Niestrawska et al. \cite{niestrawska2016} included the Cauchy stress-stretch data in both the axial and circumferential directions of 12 patch specimens. The results displayed large variability in mechanical response. However, despite the variability, when comparing the Cauchy stress at 1.15 stretch, the stresses in the circumferential direction were consistently higher than the axial direction. The experimental data chosen here was from the specimen with the behavior that was approximately median out of all the specimens (specimen AAA-1.2).

\section{Parameter Identification and Comparison of Models}
\label{chp:b5}
\subsection{Parameter identification procedure}
\label{optimization}
The parameter optimization procedure introduced by Dal et al. \cite{Dal2021} has been adopted for the constitutive models selected here for anisotropic soft tissues. The parameter identification process is conducted based on error expressions for the uniaxial tension experiments in the $\Be_1$-direction, uniaxial tension experiments in the $\Be_2$-direction, and the equibiaxial tension experiment in the $\Be_1-\Be_2$ directions,
\begin{equation*}
\begin{array}{cc}
\BcalE_{UT_1}(\Bzeta) = \sum_{i=1}^{n_{UT_1}} \left( P_{11}(\Bzeta,\lambda_i) - P_{11}^{exp}(\lambda_i)\right)^2 \quad \text{and} \quad \BcalE_{UT_2}(\Bzeta) = \sum_{i=1}^{n_{UT_2}} \left( P_{22}(\Bzeta,\lambda_i) - P_{22}^{exp}(\lambda_i)\right)^2, \\[2Ex]
\BcalE_{ET}(\Bzeta) = \sum_{i=1}^{n_{ET}} \left( P_{11}(\Bzeta,\lambda_i) - P_{11}^{exp}(\lambda_i)\right)^2 \quad \text{and} \quad \BcalE_{ET}(\Bzeta) = \sum_{i=1}^{n_{ET}} \left( P_{22}(\Bzeta,\lambda_i) - P_{22}^{exp}(\lambda_i)\right)^2, 
\end{array}
\end{equation*}
respectively, where $P_{11}$ and $P_{22}$ are the first Piola-Kirchhoff stresses, and $n_{UT_1}$, $n_{UT_2}$, $n_{ET}$ are number of data points for the UT in  $\Be_1$, UT in $\Be_2$, and ET experiments, respectively. The total cost function for the UT tests and the ET test are 
\begin{equation}
\begin{array}{rll}
\BcalE^{UT}_{TOT}(\Bzeta,\Bw) &=& w_1\BcalE_{UT_1}(\Bzeta) + w_2 \BcalE_{UT_2}(\Bzeta)\,,\\[2Ex] 
\BcalE^{ET}_{TOT}(\Bzeta,\Bw) &=& w_1\BcalE_{ET_1}(\Bzeta) + w_2 \BcalE_{ET_2}(\Bzeta) \,.
\end{array}
\label{totalcosts}
\end{equation}
The total cost functions are presented for the UT tests and ET tests individually since there is a lack of UT and ET test data that belong exactly to the same tissue. The parameters are extended to include the weights in \ref{totalcosts}, $\Bxi:= \left\{ \Bzeta,\Bw \right\}$ which is obtained from the minimization principle
\begin{equation}
\Bxi = \text{Arg}\left\{\inf_{\Bxi\in\calW} \calE_{\text{TOT}}(\Bxi)\right\} \,
\end{equation} 
along with the admissible parameter  domain
\begin{equation}
\calW = \{ \Bzeta ~\mid~ \Bzeta \in \calD ~\land~ \Bw ~\mid~   w_i \in [0,1] \}\,,
\end{equation}
where $w_1 + w_2 = 1 $. The domain $\calD$ is the physically admissible domain for the material parameters $\zeta$.
The gradient type optimization is conducted by the \verb+Fmincon+ function in Matlab where the extended cost function 
\begin{equation}
\calL(\Bxi,\Blambda^\text{eq},\Blambda^{\text{ine}})= \calE_{\text{TOT}}(\Bxi) + \sum_i\lambda_i^{\text{ine}}g^i(\Bxi) +\sum_i \lambda_i^{\text{eq}}h^i(\Bxi)
\end{equation}
is applied to equality constraints $h^i(\Bxi)$ and inequality constraints $g^i(\Bxi)$, respectively. 
To obtain an optimum solution, the variation of cost functions with respect to the extended parameters requires
\begin{equation}
\nabla_{\Bxi}\calL(\Bxi,\Blambda^\text{eq},\Blambda^{\text{ine}})=\Bnzero ~,
\end{equation}
inclusive of the Karush-Kuhn-Tucker optimality conditions for inequality constraints
\begin{equation}
\lambda_i^{\text{ine}}\ge 0
\qquad
g^i(\Bxi)\le 0
\qquad
\lambda_i^{\text{ine}}g^i(\Bxi)=0 ~,
\end{equation}
where $\lambda_i^{ine}$ are the Lagrange multipliers for the inequality constraint. The penalty parameters  $\lambda_i^{eq}$ enforce the equality constraint,
\begin{equation}
h^i(\Bxi)= 0 \,.
\end{equation}
The \verb+Fmincon+ function in Matlab
\begin{equation}
\Bxi = \textsc{Fmincon}(\calE,\Bxi_0,\BA,\Bb,\BA_{\text{eq}},\Bb_{\text{eq}})
\end{equation}
is used to minimize $\calE$, subject to the linear equality  $\BA_{\text{eq}}\Bxi=\Bb_{\text{eq}}$ and inequality $\BA\Bxi\le\Bb$. Therein, $\BA_{\text{eq}}$ is the coefficient matrix for the equality constraint, $\Bb_{\text{eq}}$ is the vector for the equality constraint, $\BA$ is the coefficient matrix for the inequality constraints, $\Bb$ is the vector for the inequality constraint, and $\Bxi_0$ are the initial points.

In this study, a hybrid optimization procedure is employed in order to reach the best parameter set. We utilize the above outlined minimization principle in conjunction with the  genetic algorithm presented in  Dal et al.~\cite{Dal2021}. The material parameter space is first trained with the genetic algorithm and the best parameter sets resulting from the genetic algorithm are used as starting points for the gradient search algorithm. 
\subsection{Comparison of hyperelastic anisotropic models}
The quality of fit metric ($\chi^2$) is used to compare the fitting performance of models. The quality of fit parameter for the uniaxial dataset of
Cooney et al.~\cite{cooney2016}, and Martins et al.~\cite{Martins2012} is
\begin{equation}\label{qofit}
\begin{array}{rl}
\chi^2 =& \sum\limits_{i=1}^{n_{UT_1}} \frac{\left(P_{11}^{UT_1}(\lambda_i) - P_{11}^{exp, UT_1}(\lambda_i)\right)^2}{P_{11}^{exp, UT_1}(\lambda_i)} + \sum\limits_{i=1}^{n_{UT_2}} \frac{\left(P_{22}^{UT_2}(\lambda_i) - P_{22}^{exp, UT_2}(\lambda_i)\right)^2}{P_{22}^{exp, UT_2}(\lambda_i)} ,
\end{array}
\end{equation}
where $n_{UT_1}$ and $n_{UT_2}$ are the number of data points and $P_{11}^{UT_1}$ and $P_{22}^{UT_2}$ are the first Piola-Kirchhoff stresses for the UT test in the $\Be_1$ and $\Be_2$ directions, respectively. Similarly, for the equibiaxial loading case, the quality of fit parameter for the fitting of the  Niestrawska et al.~\cite{niestrawska2016} dataset is

\begin{equation}\label{qofit2}
\begin{array}{rl}
\chi^2 =& \sum\limits_{i=1}^{n_{ET_1}} \frac{\left(P_{11}^{ET}(\lambda_i) - P_{11}^{exp, ET_1}(\lambda_i)\right)^2}{P_{11}^{exp, ET_1}(\lambda_i)} + \sum\limits_{i=1}^{n_{ET}} \frac{\left(P_{22}^{ET_2}(\lambda_i) - P_{22}^{exp, ET_2}(\lambda_i)\right)^2}{P_{22}^{exp, ET_2}(\lambda_i)} ,
\end{array}
\end{equation}
where $n_{ET_1}$ and $n_{ET_2}$ are the number of data points and $P_{11}^{ET_1}$ and $P_{22}^{ET_2}$ are the first Piola-Kirchhoff stresses for the ET test in the $\Be_1-\Be_2$ directions. The quality of fit metric for each model has been presented in three regions based on the stretch ranges,
g
\begin{equation}
\begin{array}{rl}
\text{region}_1 := \lambda\in[1,1/3\lambda_{max}]\,,\quad \text{region}_2 := \lambda\in[1,2/3\lambda_{max}]\,,\quad \text{region}_3 := \lambda\in[1,\lambda_{max}]\,.
\end{array}
\end{equation}

\section{Results and discussion} \label{chp:b6}
\begin{table*}[th!] 
\caption{Models sorted according to the \textsl{quality of fit} to the equibiaxial dataset of AAA tissue~\cite{niestrawska2016}, uniaxial dataset of the linea alba~\cite{cooney2016} and uniaxial dataset of the rectus sheath~\cite{Martins2012}.}\vskip -1Ex
\resizebox{1\textwidth}{!}{
\begin{tabular}{clcccp{0.25mm}lcccp{0.25mm}lccc}
\multicolumn{5}{c}{\textbf{AAA tissue}} &  & \multicolumn{4}{c}{\textbf{linea alba}}&  & \multicolumn{4}{c}{\textbf{rectus sheath}}\tabularnewline
\cline{1-5} \cline{7-10}  \cline{12-15}
\textbf{rank} & \textbf{model name} & \textbf{model type} & $\chi^{2}$ & $nop$ &  & \textbf{model name} & \textbf{model type} & $\chi^{2}$ & $nop$& & \textbf{model name} & \textbf{model type} & $\chi^{2}$ & $nop$   \tabularnewline
\cline{1-5} \cline{7-10} \cline{12-15}
1 &HNORS model~\cite{Holzapfel2015b} & GST & 2.4368 &6 & &HNORS model~\cite{Holzapfel2015b} & GST &  0.7075 &5 & &NY model~\cite{Newman1998} & $I_1, I_4$  & 0.0949 &4
\tabularnewline\cline{1-5} \cline{7-10} \cline{12-15}2 &HSGR model~\cite{Holzapfel2005} & $I_1, I_4$ & 2.4734  & 5 &  &ASMD model~\cite{Alastrue2010} & AI  &  0.9159  & 5 & &ASMD model~\cite{Alastrue2010} & AI& 0.6011 &5
\tabularnewline\cline{1-5} \cline{7-10} \cline{12-15}3 &AMDM model~\cite{Alastrue2009} & AI & 2.8541  &5 &  &AMDM model~\cite{Alastrue2009} & AI & 0.9163   &5 & &AMDM model~\cite{Alastrue2009} & AI & 0.7563 &5
\tabularnewline\cline{1-5} \cline{7-10} \cline{12-15}4 &GOH model~\cite{Gasser2006} & GST & 3.3643 &5 & &DBB model~\cite{Driessen2005} &AI & 1.0002  &6 & &DBB model~\cite{Driessen2005} & AI  & 1.8341  & 6
\tabularnewline\cline{1-5} \cline{7-10} \cline{12-15}5 &ASMD model~\cite{Alastrue2010} & AI & 3.7984  & 6 & &HSGR model~\cite{Holzapfel2005} & $I_1, I_4$  & 1.0468  &5 & &HSGR model~\cite{Holzapfel2005} & $I_1, I_4$ & 3.3545  &5
\tabularnewline\cline{1-5} \cline{7-10} \cline{12-15}6 &DBB model~\cite{Driessen2005} & AI  & 10.0814 &6 &  &OS model~\cite{Ogden2007} & $I_1, I_4$ &  1.1294 & 6 & &OS model~\cite{Ogden2007} & $I_1, I_4$& 3.3823  & 6
\tabularnewline\cline{1-5} \cline{7-10} \cline{12-15}7 &NY model~\cite{Newman1998} & $I_1, I_4$ & 11.0814  &4 &  &GOH model~\cite{Gasser2006} & GST  &  1.2228 &4 & &HNORS model~\cite{Holzapfel2015b} & GST& 3.6038   &6
\tabularnewline\cline{1-5} \cline{7-10} \cline{12-15}8 &HGO model~\cite{holzapfel2000} & $I_1, I_4$  & 47.4992  &4 &  &HGO model~\cite{holzapfel2000} & $I_1, I_4$  &  1.2529&4 & &GOH model~\cite{Gasser2006} & GST  & 3.6722  &4
\tabularnewline\cline{1-5} \cline{7-10} \cline{12-15}9 &OS model~\cite{Ogden2007} & $I_1, I_4$  & 86.1323  &5 & &NY model~\cite{Newman1998} & $I_1, I_4$ & 11.3349  &4 & &HGO model~\cite{holzapfel2000} & $I_1, I_4$  & 3.7582  &4
\tabularnewline\cline{1-5} \cline{7-10} \cline{12-15}
\end{tabular}
}
\label{sorting}
\end{table*}


The models are compared based on their quality of fit metric using the ET dataset for AAA tissue~\cite{niestrawska2016}, UT dataset for the linea alba \cite{cooney2016}, and UT dataset for the rectus sheath \cite{Martins2012}. Figures \ref{aaa-plot}--\ref{rs-plot} represent the simultaneous fit results of the models to each dataset.  In the parameter optimization procedure, we have adhered to the histological information provided by Niestrawska et al.~\cite{niestrawska2016}, Cooney et al.~\cite{cooney2016}, and Martins et al.~\cite{Martins2012} in that, for each model, we have imposed the same mean fiber directions. However, because of the lack of data on the fiber dispersion, we decided to treat the distribution parameters as model parameters to be determined  during the parameter identification procedure. In the ideal case, the distribution parameters should be obtained through fitting with histologically obtained fiber density distribution data. For nine anisotropic models, identified parameters and error bounds for AAA tissue, linea alba, and rectus sheath can be found  in Appendix.\\[1ex]
\textbf{AAA tissue:} Using the ET dataset for AAA tissue~\cite{niestrawska2016}, the models are sorted regarding the quality of fit metric, and the results are listed in Table \ref{sorting}(column \#1). The ET dataset for AAA tissue was obtained from an arterial wall specimen with two families of fibers. All models except the HGO~\cite{holzapfel2000}, and the OS~\cite{Ogden2007} models fit the AAA tissue quite successfully. The OS model fails to capture the transition of the stress-stretch relation from a linear phase to an exponential phase for the curves in both the axial and circumferential directions. On the other hand, the HGO model fail to simultaneously predict the equibiaxial stress-stretch behavior in both the circumferential and axial directions whereas they are capable of fitting either the circumferential or axial stress-stretch behavior successfully.  The first six models showed a remarkable fitting performance to the ET dataset of AAA tissue~\cite{Holzapfel2005,Alastrue2009,Alastrue2010,Newman1998,Holzapfel2015b,Gasser2006,Driessen2005}. It can be easily seen that the models taking the fiber dispersion into account are relatively more successful. For instance, HSGR and NY models use a phenomenological parameter to represent a measure of dispersion. The most successful model according to the fitting performance of the AAA tissue ET dataset is the six-parameter HNORS model \cite{Holzapfel2015b} based on the bivariate von-Mises distribution that considers the in- and  out-of-plane dispersion of fibers. The five-parameter HSGR model, and the four-parameter  NY model~\cite{Newman1998} also exhibit a remarkable fitting performance with fewer material parameters. AI-based AMDM model \cite{Alastrue2009} and GST-based GOH model \cite{Gasser2006} also demonstrated excellent fitting performance; these two models consider fiber dispersion with an identical von-Mises distribution function and use an equivalent free energy function. The only difference between these models is the integration approach. Their fitting performance is almost equivalent, however, the computational cost of the GOH model~\cite{Gasser2006} is significantly less compared to the AMDM model~\cite{Alastrue2009}. The computational cost of the  AI-based formulations are considerably higher than the GST-based ones, see also \cite{Holzapfel2017a}.
\\[1ex]
\textbf{Linea alba:} The comparison result of the models based on the UT dataset for the linea alba is listed in Table \ref{sorting}(column \#2). The UT dataset for linea alba exhibits a weak exponential form of the stress-stretch behavior. In line with the comparison for the AAA tissue dataset, fiber dispersion models have a better fitting performance for the linea alba dataset compared to invariant-based models. This is particularly due to the dispersion of fibers of the human linea alba for which the collagen fibers have a less sophisticated distribution than the collagen distribution within the arterial wall. Thus, using higher-order distribution functions to represent fiber distributions of the linea alba loses importance and increases computational complexity. For the linea alba, the HNORS model \cite{Holzapfel2015} exhibits the best fitting performance.
\\[1ex]
\textbf{Rectus sheath:} The rectus sheath tissue exhibits a relatively softer mechanical response with a weaker degree of anisotropy. The models are sorted according to the quality of fit metric and are listed in Table \ref{sorting}(column \#3). It can be observed that in both directions, the stretch-stress response of the rectus sheath tissue is nonlinear. The NY model \cite{Newman1998} has the best fitting performance for the rectus sheath dataset~\cite{Martins2012}. When considering all 9 models, the NY model \cite{Newman1998} is the only model which does not decompose the free energy function into the isotropic part and anisotropic part. Thus, even when the experimental data presents stretch-stress information in the transverse direction to the fibers, the NY model~\cite{Newman1998} considers this data as fibrous. Therefore, although the NY model \cite{Newman1998} does not represent the histological structure accurately, it is able to capture the stretch-stress response mathematically. However, the disadvantage of not decomposing the free energy function can be seen when considering the fitting performance of the NY model \cite{Newman1998} on the UT dataset for the linea alba. For the linea alba \cite{cooney2016}, the data presented in transverse to fiber direction is almost linear, and the NY model \cite{Newman1998} cannot fit this data. The next best 5 models for the rectus sheath dataset are fiber dispersion-based models. This could be due to the highly dispersed fiber structure of the rectus sheath; even transverse to the fiber direction, the experimental data shows a J-shape response. Transversely isotropic models such as the HGO model cannot capture this mechanical response, since the non-fibrous part of the tissue is modeled by the neo-Hookean model. However, another transversely isotropic model, the OS model \cite{Ogden2007} is able to show better fitting performance than the HGO model. Similar to the case of the NY model, the performance of the OS model is not due to the accurate representation of the histological structure of the tissue. The reason is that OS model uses the Gent model \cite{Gent1996} for the non-fibrous part of a tissue, which is a non-linear model. Furthermore, GST-based dispersion models have weak performance on rectus sheath data. GST-based dispersion models such as GOH~\cite{Gasser2006} and HNORS~\cite{Holzapfel2015b} models impose tension only condition for mean fiber directions. However, in the case of rectus sheath, the mean direction of fibers is along longitudinal direction and as a result GOH~\cite{Gasser2006} and HNORS~\cite{Holzapfel2015b} models are reduced to neo-Hookean model to predict mechanics behavior along transverse direction to the fibers. The same behavior can also be seen in HGO~\cite{holzapfel2000}, HSGR~\cite{Holzapfel2005} models.
\\[1ex]
\textbf{Summary:} The comparison of the results demonstrate that the fiber dispersion-based models have a superior fitting performance over the strain invariant-based models. GST-based dispersion models are more cost effective compared to the AI-based dispersion formulations regarding the computational time. On the other hand, the quality of fit of  AI- and GST-based formulations are comparable for AAA tissue, see for example the AMDM~\cite{Alastrue2009} and GOH~\cite{Gasser2006} models. However, GST-based approaches have poor performance on tissues with single family of fibers due to tension only condition. Moreover, the fitting performance of the AI-based formulations for the three distinct tissues shows that these formulations are able to adapt themselves for a variety of tissues having different families of fiber architecture and distribution.  The OS model~\cite{Ogden2007}, based on the functional form of the Gent model originally proposed for the non-Gaussian chain statistics of the rubber network, is not as successful as the models based on the Fung-type exponential representations for the anisotropic part of the free energy function, see for example HSGR~\cite{Holzapfel2005}, GOH~\cite{Gasser2006}, and AMDM~\cite{Alastrue2009} models.

\inputfig{FIGs/results_new/AAA}{aaa}
\inputfig{FIGs/results_new/LA}{la}
\inputfig{FIGs/results_new/RS}{rs}

\vskip2Ex
\textbf{Conflict of interest:} The authors declare that they have no conflict of interest.
\clearpage
\bibliographystyle{sn-chicago}
\bibliography{bibliography}%
\begin{landscape}
\section*{Appendix}
\begin{table}[th]
\caption{Identified parameters and error bounds based on AAA tissue dataset}
\resizebox{1.6\textwidth}{!}{
\begin{tabular}{llllll|lllll|lllll}
\hline
           & \multicolumn{5}{c|}{NY model}                                                          & \multicolumn{5}{c|}{HGO model}                                                        & \multicolumn{5}{c}{HSGR model}                                                      \\ \hline
Parameters & \multicolumn{5}{c|}{\begin{tabular}[c]{@{}c@{}}$k_0  =   0.1148 {~[MPa]} $ ~$ k_1  =  31.1439$  \\ $ k_2  =  1.5230e+03$~    $\varphi   =  26^{\circ}$\end{tabular}} & \multicolumn{5}{c|}{\begin{tabular}[c]{@{}c@{}}$\mu  =  2.6712 {~[MPa]}$ ~ $k_1  =  0.1742 {~[MPa]}$ ~ $k_2  =  55.9001$ \\ $\varphi   =  26^{\circ}$\end{tabular}}   & \multicolumn{5}{c}{\begin{tabular}[c]{@{}c@{}}$\mu  =  0.9347 {~[MPa]}$ ~$ k_1  =  0.2704 {~[MPa]}$ ~$ k_2  =  47.0232$ \\ $\varphi   =  26^{\circ}$ ~ $p  =  0.9126$\end{tabular}}  \\ \hline
           &                &               & \multicolumn{3}{c|}{Quality of fit}                   & \multicolumn{1}{c}{}   & \multicolumn{1}{c}{}   & \multicolumn{3}{c|}{Quality of fit} & \multicolumn{1}{c}{}   & \multicolumn{1}{c}{}  & \multicolumn{3}{c}{Quality of fit} \\ \cline{4-6} \cline{9-11} \cline{14-16} 
           & Weight         & Error         & Region 1         & Region 2         & Region 3        & Weight                 & Error                  & Region 1   & Region 2   & Region 3  & Weight                 & Error                 & Region 1   & Region 2  & Region 3  \\ \cline{2-16} 
ET-axial   & $0.1665$ & $ 0.2941 $ & $0$& $7.7801$ & $8.1602$              & $0.1$ & $ 2.4033 $ & $0$& $7.5929$ & $13.3316$         & $0.1 $ & $0.1255$ &$0$ & $1.0337$& $ 1.4823$          \\
ET-circ.    & $0.8335$ & $0.2693$ & $1.5845$ &  $0$& $2.9212$              & $0.9000$ & $14.1329$ & $8.7405$&  $0$& $34.1676$           & $0.9 $ & $0.2524 $ & $0.2511$ & $1.0337$& $0.9911 $           \\
Total      & $1.0000$ & $0.5634$ &$1.5845$ & $7.7801$& $11.0814$            & $1.0000$ & $16.5362$ & $8.7405$ & $7.5929$& $47.4992$          & $1.0000$ & $0.3779 $ & $0.2511$&$1.0337$& $ 2.4734$ \\ \hline
           & \multicolumn{5}{c|}{OS model}                                                          & \multicolumn{5}{c|}{DBB model}                                                        & \multicolumn{5}{c}{GOH model}                                                       \\ \hline
Parameters & \multicolumn{5}{c|}{\begin{tabular}[c]{@{}c@{}} $\mu  =  2.5537 {~[MPa]}$ ~$k_1  =  3.38107 {~[MPa]} $~$ J_h  =  0.1149 $\\ $J_m  =  0.2369$ ~$ \varphi  =  26^{\circ}$\end{tabular}}    & \multicolumn{5}{c|}{\begin{tabular}[c]{@{}c@{}} $G  =  2.1366 {~[MPa]} $~$ k_1  =  3.1017 {~[MPa]} $~$ k_2  =  46.8793 $\\ $\varphi   =  26^{\circ}$ ~$ \sigma   =  0.2597$ ~$ \phi_{tot}  =  0.7000$\end{tabular}}   & \multicolumn{5}{c}{\begin{tabular}[c]{@{}c@{}}$\mu  =  1.7416 {~[MPa]}$ ~$ k_1  =  4.4460 {~[MPa]} $~$ k_2  =  161.3920 {~[MPa]}$ \\ $\varphi   =  26^{\circ}$ ~$ \kappa   =  0.2256$\end{tabular}}  \\ \hline
           &                &               & \multicolumn{3}{c|}{Quality of fit}                   & \multicolumn{1}{c}{}   & \multicolumn{1}{c}{}   & \multicolumn{3}{c|}{Quality of fit} & \multicolumn{1}{c}{}   & \multicolumn{1}{c}{}  & \multicolumn{3}{c}{Quality of fit} \\ \cline{4-6} \cline{9-11} \cline{14-16} 
           & Weight         & Error         & Region 1         & Region 2         & region 3        & Weight                 & Error                  & Region 1   & Region 2   & region 3  & Weight                 & Error                 & Region 1   & Region 2  & Region 3  \\ \cline{2-16} 
ET-axial   & $0.1$ & $ 1.3169 $ & $0$& $3.7799$& $10.4829$                 & $0.1000$ & $ 1.2764 $ &$0.5860$&  $0.7356$ & $3.5558$           & $0.1000$ & $ 0.1377 $ & $0$& $1.2187$ & $1.6822$           \\
ET-circ.   & $0.9000$ & $5.7514$ & $16.3691$& $0$& $75.6493$                & $0.9000$ & $4.1076$ & $0.4600$ & $2.0767$& $7.2367$           & $0.9000$ & $0.3022$ & $0.2795$ &  $0$& $1.6821$ \\
Total      & $1.0000$ & $7.0679$ &$16.3691$ &$3.7799$& $86.1323$                & $1.0000$ & $ 5.3840$ &$1.0460$ & $2.8123$& $10.7926$         & $1.0000$ & $0.4399$ & $0.2795$ & $1.2187$& $3.3643$           \\ \hline
           & \multicolumn{5}{c|}{AMDM model}                                                        & \multicolumn{5}{c|}{ASMD model}                                                       & \multicolumn{5}{c}{HNORS model}                                                     \\ \hline
Parameters & \multicolumn{5}{c|}{\begin{tabular}[c]{@{}c@{}}$mu  =  0.9337 {~[MPa]} $~$ k_1  =  0.9118 {~[MPa]}$ ~$ k_2  =  46.8474$ \\ $\varphi   =  26^{\circ}$ ~$ b   =  3.67$\end{tabular}}    & \multicolumn{5}{c|}{\begin{tabular}[c]{@{}c@{}}$\mu  =  0.6517 {~[MPa]} $~$ k_1  =  3.5475 {~[MPa]}$ ~$ k_2  =  46.4817$ \\  $\kappa_1   =  2.3798e-07$ ~$ \kappa_2    =   0.9$ ~$ \kappa_3  =  0$\end{tabular}}   & \multicolumn{5}{c}{\begin{tabular}[c]{@{}c@{}}$\mu  =  1.8517 {~[MPa]}$ ~$ k_1  =  0.6981 {~[MPa]}$ ~$ k_2  =  59.9093$ \\ $ \kappa_{ip}   =  0.7657$~$ \kappa_{op} =  0.47$ ~$ \varphi  =  26^{\circ} $\end{tabular}}  \\ \hline
           &                &               & \multicolumn{3}{c|}{Quality of fit}                   & \multicolumn{1}{c}{}   & \multicolumn{1}{c}{}   & \multicolumn{3}{c|}{Quality of fit} & \multicolumn{1}{c}{}   & \multicolumn{1}{c}{}  & \multicolumn{3}{c}{Quality of fit} \\ \cline{4-6} \cline{9-11} \cline{14-16} 
           & Weight         & Error         & Region 1         & Region 2         & Region 3        & Weight                 & Error                  & Region 1   & Region 2   & Region 3  & Weight                 & Error                 & Region 1   & Region 2  & Region 3  \\ \cline{2-16} 
ET-axial      & $0.9000$ & $ 0.1241 $ &  $0$& $1.0682$ & $1.4882$                 & $0.9000$ & $ 0.1757 $ & $0$& $1.9907$ & $2.5057$         & $0.1000$ & $ 0.1249 $ &  $0$&  $1.0843$& $1.5269$           \\
ET-circ.   & $0.1000$ & $0.2610$ & $0.3232$ &  $0$& $1.3659$                & $0.1000$ & $1.2253$ & $0.0642$ &  $0$& $1.2907$          & $0.9000$ & $0.2429$ & $0.2309$ &  $0$& $0.9099$ \\
Total      & $1.0000$ & $0.3851$ & $0.3232$ & $1.0682$& $2.8541$                 & $1.0000$ & $1.4010$ & $0.0642$ & $1.9907$& $3.7984$          & $1.0000$ & $0.3675$ & $0.2309$ & $1.0843$& $2.4368$ \\ \hline
\end{tabular}}
\label{aaa_parameters}
\end{table}
\begin{table}[th]
\caption{Identified parameters and error bounds based on LA tissue dataset}
\resizebox{1.6\textwidth}{!}{
\begin{tabular}{llllll|lllll|lllll}
\hline
           & \multicolumn{5}{c|}{NY model}                                                          & \multicolumn{5}{c|}{HGO model}                                                        & \multicolumn{5}{c}{HSGR model}                                                      \\ \hline
Parameters & \multicolumn{5}{c|}{\begin{tabular}[c]{@{}c@{}}$k_0  =   2.9076e+03 {~[MPa]} $ ~$ k_1  =  0.0028$  \\ $ k_2  =  0.8380$~    $\varphi   =  0^{\circ}$\end{tabular}} & \multicolumn{5}{c|}{\begin{tabular}[c]{@{}c@{}}$\mu  =  1.1586 {~[MPa]}$ ~ $k_1  =  8.6064 {~[MPa]}$ ~ $k_2  =  11.4196$ \\ $\varphi   =  0^{\circ}$\end{tabular}}   & \multicolumn{5}{c}{\begin{tabular}[c]{@{}c@{}}$\mu  =  1.1525 {~[MPa]}$ ~$ k_1  =  87.4399 {~[MPa]}$ ~$ k_2  =  5.0524$ \\ $\varphi   =  0^{\circ}$ ~ $p  =  0.0458$\end{tabular}}  \\ \hline
           &                &               & \multicolumn{3}{c|}{Quality of fit}                   & \multicolumn{1}{c}{}   & \multicolumn{1}{c}{}   & \multicolumn{3}{c|}{Quality of fit} & \multicolumn{1}{c}{}   & \multicolumn{1}{c}{}  & \multicolumn{3}{c}{Quality of fit} \\ \cline{4-6} \cline{9-11} \cline{14-16} 
           & Weight         & Error         & Region 1         & Region 2         & Region 3        & Weight                 & Error                  & Region 1   & Region 2   & Region 3  & Weight                 & Error                 & Region 1   & Region 2  & Region 3  \\ \cline{2-16} 
ET-axial   & $0.1000$ & $0.0399 $ & $0.6616$& $1.5653$ & $1.7889$              & $0.1000$ & $ 0.0079 $ &$0.3335$& $0.4445$& $0.4892$         & $0.1001 $ & $0.0079$ &$0.3337$ & $0.4448$& $ 0.4895$          \\
ET-circ.    & $0.9000$ & $0.6957$ & $2.6033$ &  $8.3558$& $9.5460$              & $0.9000$ & $0.0171$ & $0.7246$ &$0.7392$& $0.7637$           & $0.8999 $ & $0.0123 $ & $0.5304$ & $0.5411$& $0.5573 $           \\
Total      & $1.0000$ & $0.7356$ &$3.2648$ & $9.9211$& $11.3349$            & $1.0000$ & $0.0250$ &$1.0581$&$1.1837$& $1.2529$          & $1.0000$ & $0.0199 $ & $0.8641$&$0.9859$& $ 1.0468$ \\ \hline
           & \multicolumn{5}{c|}{OS model}                                                          & \multicolumn{5}{c|}{DBB model}                                                        & \multicolumn{5}{c}{GOH model}                                                       \\ \hline
Parameters & \multicolumn{5}{c|}{\begin{tabular}[c]{@{}c@{}} $\mu  =  1.7592 {~[MPa]}$ ~$k_1  =  9.8814 {~[MPa]} $~$ J_h  =  0.2883 $\\ $J_m  =  0.0944$ ~$ \varphi  =  0^{\circ}$\end{tabular}}    & \multicolumn{5}{c|}{\begin{tabular}[c]{@{}c@{}} $G  =  0.0500 {~[MPa]} $~$ k_1  =  56.0009 {~[MPa]} $~$ k_2  =  0.7921 $\\ $\varphi   =  0^{\circ}$ ~$ \sigma   =  0.4985$ ~$ \phi_{tot}  =  0.5200$\end{tabular}}   & \multicolumn{5}{c}{\begin{tabular}[c]{@{}c@{}}$\mu  =  2.3048 {~[MPa]}$ ~$ k_1  =  17.2552 {~[MPa]} $~$ k_2  =  20.5490 {~[MPa]}$ \\ $\varphi   =  0^{\circ}$ ~$ \kappa   =  0.0998$\end{tabular}}  \\ \hline
           &                &               & \multicolumn{3}{c|}{Quality of fit}                   & \multicolumn{1}{c}{}   & \multicolumn{1}{c}{}   & \multicolumn{3}{c|}{Quality of fit} & \multicolumn{1}{c}{}   & \multicolumn{1}{c}{}  & \multicolumn{3}{c}{Quality of fit} \\ \cline{4-6} \cline{9-11} \cline{14-16} 
           & Weight         & Error         & Region 1         & Region 2         & region 3        & Weight                 & Error                  & Region 1   & Region 2   & region 3  & Weight                 & Error                 & Region 1   & Region 2  & Region 3  \\ \cline{2-16} 
UT-trans.   & $0.1000$ 	& $ 8.4265e-04 $ & $0.0709$& $0.0736$& $0.0797$                 & $0.0100$ & $ 0.0030 $ &$0.0140$&  $0.0754$ & $0.1022$           & $0.1000$ & $ 0.0076 $ & $0.0335$& $0.4444$ & $0.4892$           \\
ET-circ.   & $0.9000$ & $0.0245$ & $0.9835$& $1.0114$& $1.0496$                & $0.9900$ & $0.0199$ & $0.8302$ & $0.8902$& $0.8980$           & $0.9000$ & $0.0164$ & $0.6971$ &  $0.7106$& $0.7336$ \\
Total      & $1.0000$ & $0.0253$ &$1.0544$ &$1.0851$& $1.1294$                & $1.0000$ & $ 0.0299$ &$0.8441$ & $0.9656$& $1.0002$         & $1.0000$ & $0.024$ & $1.0306$ & $1.1551$& $1.2228$           \\ \hline
           & \multicolumn{5}{c|}{AMDM model}                                                        & \multicolumn{5}{c|}{ASMD model}                                                       & \multicolumn{5}{c}{HNORS model}                                                     \\ \hline
Parameters & \multicolumn{5}{c|}{\begin{tabular}[c]{@{}c@{}}$mu  =  0.7418 {~[MPa]} $~$ k_1  =  15.3359 {~[MPa]}$ ~$ k_2  =  10.6226$ \\ $\varphi   =  0^{\circ}$ ~$ b   =  2.6374$\end{tabular}}    & \multicolumn{5}{c|}{\begin{tabular}[c]{@{}c@{}}$\mu  =  0.7420 {~[MPa]} $~$ k_1  =  15.6560 {~[MPa]}$ ~$ k_2  =  10.6286$ \\  $\kappa_1   =  4.999e-07$ ~$ \kappa_2    =   5.2757$ ~$ \kappa_3  = 0.0000$\end{tabular}}   & \multicolumn{5}{c}{\begin{tabular}[c]{@{}c@{}}$\mu  =  2.3050 {~[MPa]}$ ~$ k_1  =  270.6246 {~[MPa]}$ ~$ k_2  =  5.0500$ \\ $ \kappa_{ip}   =  0.6448$~$ \kappa_{op} =  0.4999$ ~$ \varphi  =  0^{\circ} $\end{tabular}}  \\ \hline
           &                &               & \multicolumn{3}{c|}{Quality of fit}                   & \multicolumn{1}{c}{}   & \multicolumn{1}{c}{}   & \multicolumn{3}{c|}{Quality of fit} & \multicolumn{1}{c}{}   & \multicolumn{1}{c}{}  & \multicolumn{3}{c}{Quality of fit} \\ \cline{4-6} \cline{9-11} \cline{14-16} 
           & Weight         & Error         & Region 1         & Region 2         & Region 3        & Weight                 & Error                  & Region 1   & Region 2   & Region 3  & Weight                 & Error                 & Region 1   & Region 2  & Region 3  \\ \cline{2-16} 
Axial      & $0.1000$ & $ 0.0012 $ &$0.0800$& $0.0843$& $0.0944$                & $0.1000$ & $ 0.0012 $ &$0.0800$&$0.0844$& $0.0944$         & $0.8998$ & $ 0.0076 $ &  $0.3255$&  $0.4304$& $0.4789$           \\
ET-circ.   & $0.9000$ & $0.0170$ & $0.7877$&$0.8052$& $0.8219$                & $0.9000$ & $0.0170$ &$0.7873$&$0.0847$& $0.8214$          & $0.1002$ & $0.0055$ & $0.2138$ &  $0.2221$& $0.2286$ \\
Total      & $1.0000$ & $0.0182$ &$0.8677$&$0.8895$& $0.9163$                 & $1.0000$ & $0.0182$ &$0.8673$&$0.8890$& $0.9159$          & $1.0000$ & $0.0131$ & $0.5494$ & $0.6525$& $0.7075$ \\ \hline
\end{tabular}}
\label{la_parameters}
\end{table}
\begin{table}[th]
\caption{Identified parameters and error bounds based on RS tissue dataset}
\resizebox{1.6\textwidth}{!}{
\begin{tabular}{llllll|lllll|lllll}
\hline
           & \multicolumn{5}{c|}{NY model}                                                          & \multicolumn{5}{c|}{HGO model}                                                        & \multicolumn{5}{c}{HSGR model}                                                      \\ \hline
Parameters & \multicolumn{5}{c|}{\begin{tabular}[c]{@{}c@{}}$k_0  =  0.8366 {~[MPa]} $ ~$ k_1  =  6.3492$  \\ $ k_2  =  25.7852$~    $\varphi   =  90^{\circ}$\end{tabular}} & \multicolumn{5}{c|}{\begin{tabular}[c]{@{}c@{}}$\mu  =  0.0700 {~[MPa]}$ ~ $k_1  =  0.4465 {~[MPa]}$ ~ $k_2  =  7.6186$ \\ $\varphi   =  90^{\circ}$\end{tabular}}   & \multicolumn{5}{c}{\begin{tabular}[c]{@{}c@{}}$\mu  =  0.0900 {~[MPa]}$ ~$ k_1  =  3.1260 {~[MPa]}$ ~$ k_2  =  35.2661$ \\ $\varphi   =  90^{\circ}$ ~ $p  =  0.1000$\end{tabular}}  \\ \hline
           &                &               & \multicolumn{3}{c|}{Quality of fit}                   & \multicolumn{1}{c}{}   & \multicolumn{1}{c}{}   & \multicolumn{3}{c|}{Quality of fit} & \multicolumn{1}{c}{}   & \multicolumn{1}{c}{}  & \multicolumn{3}{c}{Quality of fit} \\ \cline{4-6} \cline{9-11} \cline{14-16} 
           & Weight         & Error         & Region 1         & Region 2         & Region 3        & Weight                 & Error                  & Region 1   & Region 2   & Region 3  & Weight                 & Error                 & Region 1   & Region 2  & Region 3  \\ \cline{2-16} 
UT-longi.   & $0.1015$ & $2.7021e-04 $ & $0.0252$& $0.0279$ & $0.0287$              & $0.1000$ & $ 0.0059 $ & $0.9175$& $1.0585$ & $1.0908$         & $0.1000 $ & $0.0048$ &$0.7407$ & $0.8783$& $ 0.9004$          \\
UT-trans.    & $0.8985$ & $6.8242e-04$ & $0.0275$ &  $0.0557$& $0.0663$              & $0.9000$ & $0.0981$ & $0.0122$&  $0.4693$& $2.6674$           & $0.9000 $ & $0.0863 $ & $0.0059$ & $0.3937$& $2.4541 $           \\
Total      & $1.0000$ & $9.5263e-04$ &$0.0526$ & $0.0836$ & $0.0949$            & $1.0000$ & $0.0981$ & $0.9297$ & $1.5278$& $3.7582$          & $1.0000$ & $0.0911 $ & $0.7467$&$1.2720$& $ 3.3545$ \\ \hline
           & \multicolumn{5}{c|}{OS model}                                                          & \multicolumn{5}{c|}{DBB model}                                                        & \multicolumn{5}{c}{GOH model}                                                       \\ \hline
Parameters & \multicolumn{5}{c|}{\begin{tabular}[c]{@{}c@{}} $\mu  =  0.6971 {~[MPa]}$ ~$k_1  =  0.3315 {~[MPa]} $~$ J_h  =  0.1381 $\\ $J_m  =  0.3263$ ~$ \varphi  =  90^{\circ}$\end{tabular}}    & \multicolumn{5}{c|}{\begin{tabular}[c]{@{}c@{}} $G  =  0.05 {~[MPa]} $~$ k_1  =  9.0099 {~[MPa]} $~$ k_2  =  5.0009 $\\ $\varphi   =  90^{\circ}$ ~$ \sigma   =  1.0473$ ~$ \phi_{tot}  =  0.3001$\end{tabular}}   & \multicolumn{5}{c}{\begin{tabular}[c]{@{}c@{}}$\mu  =  0.0700 {~[MPa]}$ ~$ k_1  =  23.5547 {~[MPa]} $~$ k_2  =  131.2277 {~[MPa]}$ \\ $\varphi   =  90^{\circ}$ ~$ \kappa   =  0.2500$\end{tabular}}  \\ \hline
           &                &               & \multicolumn{3}{c|}{Quality of fit}                   & \multicolumn{1}{c}{}   & \multicolumn{1}{c}{}   & \multicolumn{3}{c|}{Quality of fit} & \multicolumn{1}{c}{}   & \multicolumn{1}{c}{}  & \multicolumn{3}{c}{Quality of fit} \\ \cline{4-6} \cline{9-11} \cline{14-16} 
           & Weight         & Error         & Region 1         & Region 2         & region 3        & Weight                 & Error                  & Region 1   & Region 2   & region 3  & Weight                 & Error                 & Region 1   & Region 2  & Region 3  \\ \cline{2-16} 
ET-axial   & $0.9000$ 	& $ 0.0173 $ & $2.2566$& $2.6379$& $2.7472$                 & $0.9900$ & $ 0.0082 $ &$1.2097$&  $1.54155$ & $1.5579$           & $0.0100$ 	& $ 0.0037 $ & $0.4653$& $0.5722$& $0.5904$           \\
ET-circ.   & $0.1000$ & $0.0021$ & $0.4551$& $0.6240$& $0.6351$               & $0.0100$ & $0.0059$ & $0.1067$ & $0.1365$& $0.2762$                 & $0.9900$ & $0.1033$ & $0.0365$& $0.6314$& $3.0817$ \\
Total      & $1.0000$ & $0.0194$ &$2.7117$ &$3.2618$& $3.3823$                & $1.0000$ & $ 0.0141$ &$1.3164$ & $1.6780$& $1.8341$            & $1.0000$ & $0.1059$ &$0.5019$ &$1.2036$& $3.6722$            \\ \hline
           & \multicolumn{5}{c|}{AMDM model}                                                        & \multicolumn{5}{c|}{ASMD model}                                                       & \multicolumn{5}{c}{HNORS model}                                                     \\ \hline
Parameters & \multicolumn{5}{c|}{\begin{tabular}[c]{@{}c@{}}$mu  =  1.8478e-05 {~[MPa]} $~$ k_1  =  3.1303 {~[MPa]}$ ~$ k_2  =  9.3118$ \\ $\varphi   =  90^{\circ}$ ~$ b   =  0.2620$\end{tabular}}    & \multicolumn{5}{c|}{\begin{tabular}[c]{@{}c@{}}$\mu  =  0.0503 {~[MPa]} $~$ k_1  =  1.9067{~[MPa]}$ ~$ k_2  =  13.5294$ \\  $\kappa_1   =  6.3935e-07$ ~$ \kappa_2    =   8.9346e-07$ ~$ \kappa_3  =  0$\end{tabular}}   & \multicolumn{5}{c}{\begin{tabular}[c]{@{}c@{}}$\mu  = 0.0851 {~[MPa]}$ ~$ k_1  =  16.7116 {~[MPa]}$ ~$ k_2  =  152.3987$ \\ $ \kappa_{ip}   =  0.5217$~$ \kappa_{op} =  0.3900$ ~$ \varphi  =  90^{\circ} $\end{tabular}}  \\ \hline
           &                &               & \multicolumn{3}{c|}{Quality of fit}                   & \multicolumn{1}{c}{}   & \multicolumn{1}{c}{}   & \multicolumn{3}{c|}{Quality of fit} & \multicolumn{1}{c}{}   & \multicolumn{1}{c}{}  & \multicolumn{3}{c}{Quality of fit} \\ \cline{4-6} \cline{9-11} \cline{14-16} 
           & Weight         & Error         & Region 1         & Region 2         & Region 3        & Weight                 & Error                  & Region 1   & Region 2   & Region 3  & Weight                 & Error                 & Region 1   & Region 2  & Region 3  \\ \cline{2-16} 
Axial      & $0.9000$ & $ 0.0065 $ &  $.4221$& $0.4869$ & $0.5097$                 & $0.9000$ & $ 0.0347 $ & $0.0757$& $0.1804$ & $0.4085$           & $0.8999$ & $ 0.0051 $ &  $0.4923$&  $0.5730$& $0.6161$           \\
ET-circ.   & $0.1000$ & $9.1167e-04$ & $0.1763$ &  $0.2386$& $0.2466$                & $0.1000$ & $0.0024$ & $0.1200$ &  $0.1478$& $0.1925$          & $0.1001$ & $0.1008$ & $0.0299$ &  $0.5931$& $2.9878$ \\
Total      & $1.0000$ & $7.41167e-03$ & $0.5984$ & $0.7254$& $0.7563$                 & $1.0000$ & $0.0371$ & $0.1957$ & $0.3282$& $0.6011$          & $1.0000$ & $0.1059$ & $0.5221$ & $1.1661$& $3.6038$ \\ \hline
\end{tabular}}
\label{aaa_parameters}
\end{table}

\end{landscape}
\end{document}